\documentclass[11pt]{article}

\usepackage[T1,T2A]{fontenc}
\usepackage[utf8]{inputenc}
\usepackage{amsfonts}
\usepackage{amssymb}
\usepackage{amsmath}
\usepackage{epsfig}
\usepackage{color}  
\usepackage{theorem}
\usepackage{tikz}

\newtheorem{theorem}{Theorem}[section]
\newtheorem{lemma}[theorem]{Lemma}
\newtheorem{corollary}[theorem]{Corollary}
\newtheorem{proposition}[theorem]{Proposition}

{\theorembodyfont{\rmfamily} 

}

\usepackage[color links=false, backref=page]{hyperref} 

\def\Q{\mathbb{Q}}

\def\F{\mathbb{F}}

\def\H{\mathbb{H}}

\def\M{\mathcal{M}}
\def\N{\mathcal{N}}
\def\S{\mathcal{S}}
\def\U{\mathcal{U}}

\def\PSL{\rm{PSL}}
\def\A{\rm{Aut}}

\title{Regularity properties of Macbeath--Hurwitz and related maps and surfaces}
\date{\today}
\author{Gareth A. Jones}

\begin{document}

\medskip

\maketitle

\medskip

\begin{abstract}
The Macbeath--Hurwitz maps $\M$ of type $\{3,7\}$, obtained from the Hurwitz groups $G={\rm PSL}_2(q)$ found by Macbeath, are fully regular by a result of Singerman, with automorphism group $G\times{\rm C}_2$ or ${\rm PGL}_2(q)$. Hall's criterion determines which of these two properties, called inner and outer regularity, $\M$ has. Inner (but not outer) regular maps $\M$ yield non-orientable regular maps $\M/{\rm C}_2$ of the same type with automorphism group $G$.
If $q=p^3$ for a prime $p\equiv\pm 2$ or $\pm 3$ mod~$(7)$ the unique map $\M$ is inner regular if and only if $p\equiv 1$ mod~$(4)$. If $q=p$ for a prime $p\equiv\pm 1$ mod~$(7)$ there are three maps $\M$; we use the density theorems of Frobenius and Chebotarev to show that in this case the sets of such primes $p$ for which $0, 1, 2$ or $3$ of them are inner regular have relative densities $1/8$, $3/8$, $3/8$ and $1/8$ respectively. Hall's criterion and its consequences are extended to the analogous Macbeath maps of type $\{3,n\}$ obtained from ${\rm PSL}_2(q)$ for all $n\ge 7$; theoretical predictions on their number and properties are supported by evidence from the map databases of Conder and Poto\v cnik.
\end{abstract}

  
\medskip



\section{Introduction}

 Hurwitz~\cite{Hur93} showed that if $\S$ is a compact Riemann surface (equivalently, a complex projective algebraic curve) of genus $g\ge 2$ then its automorphism group $G={\rm Aut}\,\S$ has order $|G|\le 84(g-1)$. The surfaces $\S$ and groups $G$ attaining this bound are called {\sl Hurwitz surfaces\/} and {\sl Hurwitz groups}.  For a long time the only known examples were Klein's quartic curve~\cite{Kle78} of genus~$3$ with automorphism group ${\rm PSL}_2(7)$, and an example of genus~$7$ with automorphism group ${\rm PSL}_2(8)$ discovered by Fricke~\cite{Fri}; this example was long forgotten, but was rediscovered in 1965 by Macbeath~\cite{Macb65}. In 1961 he showed in~\cite{Macb61} that any Hurwitz surface yields infinitely many more as its mod~$(m)$ homology covers for integers $m\ge 2$ (the `Macbeath trick'). In 1969 Macbeath~\cite{Macb69} extended the examples due to Klein and Fricke by proving the following:

\begin{theorem}[Macbeath]\label{th:Macb}
The group $G={\rm PSL}_2(q)$ is a Hurwitz group if and only if
\begin{enumerate}
\item $q=7$, with a unique Hurwitz surface $\S$, or
\item $q=p$ for some prime $p\equiv\pm 1$ {\rm mod}~$(7)$, with three Hurwitz surfaces $\S$  for each such $q$, or
\item $q=p^3$ for some prime $p\equiv\pm 2$ or $\pm 3$ {\rm mod}~$(7)$, with one Hurwitz surface $\S$ for each such $q$.
\end{enumerate}
\end{theorem}

Any Hurwitz surface $\S$ carries an orientably regular map $\M$ of type $\{3,7\}$ (i.e.~a 7-valent triangulation), called a {\sl Hurwitz map}, with orientation-preserving automorphism group ${\rm Aut}^+\M=G={\rm Aut}\,\S$.
The maps $\M$ on the surfaces in Theorem~\ref{th:Macb} are called {\sl Macbeath--Hurwitz maps},. It follows from a result of Singerman~\cite{Sin74} that they are fully regular, with automorphism group
\[{\rm Aut}\,\M\cong{\rm PSL}_2(q)\times{\rm C}_2\quad\hbox{or}\quad{\rm PGL}_2(q);\]
we will say that $\M$ and $\S$ are {\sl inner} or {\sl outer regular} respectively. In either case, the orientation-reversing automorphisms of $\M$ act as anticonformal automorphisms of $\S$. When $\S$ is inner regular, its central quotient $\S/{\rm C}_2$ is a non-orientable Klein surface of genus $g'=g+1$ with automorphism group $G$ attaining the upper bound of order $84(g'-2)$ for such a surface of this genus (such groups are called $H^*$-{\sl groups}), while $\M/{\rm C}_2$ is a non-orientable regular map of type $\{3,7\}$ with automorphism group $G$ and with orientable double cover $\M$.

 It is therefore of interest to know whether a given Macbeath--Hurwitz map $\M$ (or surface $\S$) is inner or outer regular. For example, the map in case~(1) of Theorem~\ref{th:Macb} (Klein's map of genus $3$) is outer regular, whereas the Macbeath--Fricke map of genus $7$ (case~(3) with $p=2$) is inner regular. Hall~\cite[Theorem~2.9]{Hall77} proved the following (restated here in terms of maps, rather than Riemann and Klein surfaces):
 
 \begin{theorem}[Hall's Criterion]\label{th:Hall}
 Let $\M$ be a Macbeath--Hurwitz map with orientation-preserving automorphism group ${\rm PSL}_2(q)$, and let $\pm t$ be the trace of an automorphism rotating $\M$ by $2\pi/7$ around a vertex. Then $\M$ is inner or outer regular as $3-t^2$ is or is not a square in ${\mathbb F}_q$.
 \end{theorem}
 
 She used this to show that the maps in case~(3) of Theorem~\ref{th:Macb} are inner regular if and only if $p=2$ or $p\equiv 1$ mod~$(4)$. However, the situation is more complicated in case~(2): for example, if $p=167$, $13$, $43$ or $181$ then the number of inner regular maps is $0$, $1$, $2$ or $3$ respectively. (The fact that there can be both inner and outer regular maps for the same prime $p$ shows that these properties are not  Galois-invariant: Streit~\cite{Str} has shown that for each prime $p$ in case~(2) the three curves $\S$ are defined over the real cyclotomic field $\Q(\cos 2\pi/7)$ and are conjugate under the Galois group ${\rm C}_3$ of that field.)

 By applying the density theorems of Frobenius and Chebotarev to Hall's criterion we will prove the following, our first main theorem. For each prime $p\equiv \pm 1$ mod~$(7)$, we define $k_p$ to be the number of inner regular Macbeath--Hurwitz maps with orientation-preserving automorphism group ${\rm PSL}_2(p)$:
 
 \begin{theorem}[Density Theorem for Macbeath--Hurwitz maps]\label{th:densities}
The sets $\Sigma_k$ of primes $p\equiv\pm 1$ {\rm mod}~$(7)$ for which $k_p=k=0$, $1$, $2$ or $3$ have relative densities $1/8$, $3/8$, $3/8$ and $1/8$ respectively.
 \end{theorem}
 
This is supported by computational data from Marston Conder for the first $400$ primes $p\equiv\pm 1$ mod~$(7)$, where the proportions of primes with $k_p=0$, $1$, $2$ or $3$ are close to these values (see Section~\ref{sec:regularity} and the Appendix). It would be good to have a complete characterisation of the primes achieving each of these values of $k_p$, but at present the only general result we have is Hall's corollary to her criterion~\cite{Hall77}:

\begin{theorem}[Hall's Parity Theorem]\label{th:parity}
$k_p$ is odd or even as $p\equiv\pm 1$ {\rm mod}~$(4)$.
\end{theorem}

This has the following immediate consequences:

\begin{corollary}
\begin{enumerate}
\item If $p\equiv 1$ or $13$ {\rm mod}~$(28)$ there is at least one inner regular Macbeath--Hurwitz map $\M$ with ${\rm Aut}^+\M\cong{\rm PSL}_2(p)$;
\item If $p\equiv -1$ or $-13$ {\rm mod}~$(28)$ there is at least one outer regular Macbeath--Hurwitz map $\M$ with ${\rm Aut}^+\M\cong{\rm PSL}_2(p)$;
\item If $p\equiv 5$, $9$, $17$ or $25$ {\rm mod}~$(28)$ the unique Macbeath--Hurwitz map $\M$ with ${\rm Aut}^+\M\cong{\rm PSL}_2(p^3)$ is inner regular;
\item If $p\equiv -5$ $-9$, $-17$ or $-25$ {\rm mod}~$(28)$ the unique Macbeath--Hurwitz map $\M$ with ${\rm Aut}^+\M\cong{\rm PSL}_2(p^3)$ is outer regular.
\end{enumerate}
\end{corollary}

In the second half of this paper we find analogues of these results for a more extensive class of triangular maps. Hall proved her regularity criterion 
only for the Macbeath--Hurwitz maps, but in fact her proof applies, with only minor modifications, to a much wider family of maps of type $\{3,n\}$ for all $n\ge 7$. Although her criterion has subsequently been generalised even further by Conder, Poto\v cnik and \v Sir\'a\v n~\cite{CPS}, we will, largely for historical reasons, give a short outline of her unpublished proof here, in Section~\ref{sec:Hall}, presented to show how it applies to all $n\ge 7$.

Let us define $\M$ to be a {\sl Macbeath map of type} $\{3,n\}$, or simply a {\sl Macbeath map} if the type is understood, if it is an orientably regular map of that type with orientation-preserving automorphism group $G={\rm PSL}_2(q)$ for some prime power $q$ coprime to $n$, where $n\ge 7$.
By Singerman's result~\cite{Sin74}  $\M$ is fully regular, and as an extension of Theorem~\ref{th:Hall} we have the following:

 \begin{theorem}[Extended Regularity Criterion]\label{th:extcriterion}
 Let $\M$ be a Macbeath map of type $\{3,n\}$ with orientation-preserving automorphism group ${\rm PSL}_2(q)$, and let $\pm t$ be the trace of an automorphism rotating $\M$ by $2\pi/n$ around a vertex. Then $\M$ is inner or outer regular as $3-t^2$ is or is not a square in ${\mathbb F}_q$.
 \end{theorem}

These Macbeath maps exist if and only if $q=p^d$ where the prime $p$ has multiplicative order $d$ in ${\mathbb Z}_n^*/\langle\pm 1\rangle$. Hall's parity theorem has a natural, but slightly more complicated, generalisation to them. Let $\Psi_n$ be the minimal polynomial, over $\mathbb Q$, of the algebraic integer $2\cos(2\pi/n)$.

\begin{theorem}[Extended Parity Theorem]\label{th:extparity}
 Let $G={\rm PSL}_2(q)$ where $q=p^d$ for some prime $p$ and odd $d$. Then the number $l_p$ of outer regular Macbeath maps of type $\{3,n\}$ with orientation-preserving automorphism group $G$ is even or odd as $\Psi_n(1)$ is or is not a square in ${\mathbb F}_q$.
\end{theorem}

For instance, since $\Psi_7(1)=-1$, which is a square in ${\mathbb F}_p$ if and only if $p\not\equiv -1$ mod~$(4)$, and since $k_p$ and $l_p$ have opposite parities when $n=7$, we obtain Theorem~\ref{th:parity} as a special case. Note that Theorem~\ref{th:extparity} does not apply when $d$ is even; however, it does apply in the  most frequently occurring case, when $d=1$, and more generally when $\varphi(n)/2$ is odd. This result is proved in Section~\ref{sec:parity}.

Sections~\ref{sec:oddexamples} and \ref{sec:evenexamples} contain illustrative applications of Theorems~\ref{th:extcriterion} and \ref{th:extparity} to small odd and even values of $n$, with the latter treated separately on account of some extra technical complications arising. The theoretical predictions given by these theorems are compared, for genera up to 1501, with the evidence in the extensive databases of regular maps~\cite{Con09} and \cite{Pot} compiled by Conder and by Conder and Poto\v cnik. In all cases examined, they are consistent with this evidence.

Applying the density theorems mentioned earlier requires knowledge of the Galois groups of certain polynomials; these are discussed in more detail in Section~\ref{sec:Galois}. In Section~\ref{sec:m>3} we outline how the methods used here can be extended to Macbeath maps of type $\{m,n\}$ for any fixed $m>3$, using the case $m=4$ as an example. Some further generalisations of this work are briefly discussed in Section~\ref{sec:general}, and Section~\ref{sec:Conder} contains Conder's computational data on Macbeath--Hurwitz maps and groups.

Although the main objects of study in this paper, namely certain families of highly symmetric Riemann and Klein surfaces and their triangulations, are geometric in nature, the fact that their automorphism groups are most conveniently defined over finite fields means that the algebraic and number-theoretic properties of these fields and of some related polynomials play a major role in this work. In fact, one of the main attractions of this topic is the way in which it brings together ideas from different areas of mathematics. For example, the fact that the Riemann surfaces considered here are all uniformised by subgroups of finite index in triangle groups means that they and the maps they carry are all instances of Grothendieck's dessins d'enfants (see~\cite{JW, JZ}, for example), with orientable regularity of the maps equivalent to regularity of the dessins.
 
\bigskip
 
\centerline{\bf \large Acknowledgments}

\medskip

The author is very grateful to Marston Conder for computational data about regular maps, some presented in the Appendix, to Martin Ma\v caj for computational results about Galois groups, and to Jozef \v Sir\'a\v n and Primo\v z Poto\v cnik for helpful comments.

 
 \section{Triangle groups and regular maps}\label{sec:triangle}
 
 If $\frac{1}{m}+\frac{1}{n}<\frac{1}{2}$, let $T$ be a triangle $ABC$ in the hyperbolic plane $\H$ with internal angles $\pi/n$, $\pi/2$ and $\pi/m$ at $A$, $B$ and $C$. Then the extended triangle group $\Delta=\Delta[2,m,n]$ of type $(2,m,n)$ is the group generated by reflections $a, b$ and $c$ of $\H$ in the sides of $T$ opposite $A$, $B$ and $C$. It follows from the fact that $\H$ is simply connected that $\Delta$ has a presentation
 \[\Delta=\langle a, b, c\mid a^2=b^2=c^2=(ab)^m=(bc)^n=(ca)^2=1\rangle.\]
The ordinary triangle group $\Delta^+=\Delta(2,m,n)$ is the orientation-preserving subgroup of index $2$ in $\Delta$ consisting of the elements of even word-length in the generators $a$, $b$ and $c$. It has a presentation
 \[\Delta^+=\langle X, Y, Z\mid X^m=Y^n=Z^2=XYZ=1\rangle\]
 where $X=ab$, $Y=bc$ and $Z=(XY)^{-1}=ca$ are rotations through $2\pi/m$, $2\pi/n$ and $\pi$ around $C$, $A$ and $B$. The images of the side $AB$, under either of these two groups, form the edges of a map $\U=\U(m,n)$ of type $\{m,n\}$ on $\H$, which is invariant under both groups: the vertices of $\U$ (the images of $A$) have valency $n$, the edges are hyperbolic geodesics, and the faces (each with an image of $C$ at its midpoint) are hyperbolic $m$-gons. The tessellation of $\H$ by the images of $T$ is the barycentric subdivision of $\U$.
 
 If $K$ is a torsion-free normal subgroup of $\Delta$ then $\U/K$ is a map $\M$ of type $\{m,n\}$ on the surface $\H/K$. This map $\M$ is regular, in the sense that its automorphism group $\A\,\M$, isomorphic to $\Delta/K$, acts transitively on its vertex-edge-face flags. If $K\le\Delta^+$ then $\H/K$ is a Riemann surface and $\M$ is orientable, whereas if $K\not\le\Delta^+$ then $\H/K$ is a Klein surface without boundary, and $\M$ is non-orientable. If we merely assume that $K$ is normal in $\Delta^+$, rather than $\Delta$, then $\M$ is orientably regular, meaning that its orientation-preserving automorphism group $\A^+\M$, isomorphic to $\Delta^+/K$, acts transitively on its arcs (directed edges); then $\M$ is (fully) regular if $K$ is normal in $\Delta$, but if not then $K$ is one of a pair of normal subgroups of $\Delta^+$, transposed by conjugation in $\Delta$ and corresponding to a chiral (mirror-image) pair of orientably regular maps.
 
Conversely, if $\M$ is any regular map $\M$ of type $\{m,n\}$ then $\A\,\M$ is generated by three reflections $a'$, $b'$ and $c'$ changing the $i$-dimensional component of a specific flag of $\M$ for $i=0$, $1$ and $2$, while preserving the other two components. These generators satisfy all the defining relations of $\Delta$, and their products have the same orders, so there is an epimorphism $\Delta\to\A\,\M$ with a torsion-free kernel $K$, and $\M\cong\U/K$ with $\A\,\M\cong\Delta/K$. Thus one can regard $\U$ as the universal map of type $\{m,n\}$, in the sense that every regular map of this type is isomorphic to a quotient of it.
 (In fact, this applies to {\em every} map of this type, regular or not, since every map is a quotient of a regular map of the same type, by some group of automorphisms; however, we will not need this full generality here.) The underlying surface $\H/K$ of $\M$ is compact if and only if $K$ has finite index in $\Delta$, as we will assume throughout this paper.
 
 There is an obvious analogue of the preceding paragraph for orientably regular maps, with the orientation-preserving automorphism group $\A^+\M$ replacing $\A\,\M$, and $\Delta^+$ replacing $\Delta$.

All of the statements in this section extend to the case where $\frac{1}{m}+\frac{1}{n}\ge\frac{1}{2}$, provided one replaces $\H$ with the complex plane in the case of equality, or the Riemann sphere in the case of strict inequality. However, the resulting Riemann and Klein surfaces, together with the regular maps on them, are so well known and understood that we will avoid them here.
Similarly, one can extend these statements to arbitrary cocompact triangle groups $\Delta$, without the restriction that $(ca)^2=1$, provided one replaces maps with hypermaps, surface embeddings of hypergraphs. Here we will restrict our attention to maps, though extending some of our results to hypermaps would be an interesting project.
 

\section{Hurwitz surfaces, groups and maps}\label{sec:Hurwitz}

Let $\S$ be a compact Riemann surface of genus $g\ge 2$. Hurwitz's upper bound $|G|\le 84(g-1)$ for $G={\rm Aut}\,\S$ is attained if and only if $\S\cong \H/K$, where $\H$ is the hyperbolic plane and $K$ is a proper normal subgroup of finite index in the (orientation-preserving) triangle group
\[\Delta^+=\Delta(2,3,7)=\langle X, Y \mid X^3=Y^7=Z^2=XYZ=1\rangle,\]
in which case $G\cong \Delta/K$. Since $2$, $3$ and $7$ are distinct primes, such a subgroup $K$ is automatically torsion-free, and is therefore a surface group of genus $g$.
It follows that a non-trivial finite group $G$ is a Hurwitz group (that is, it attains the Hurwitz bound for some surface $\S$) if and only if it has generators $x$, $y$ and $z$ (the images of $X$, $Y$ and $Z$) satisfying
\[x^3=y^7=z^2=xyz=1.\]
(See~\cite{Con90, Con10} for comprehensive surveys on Hurwitz groups.)

For any such group $G$ there is a corresponding Hurwitz map $\M=\U/K$ of type $\{3,7\}$ on $\S$, where $\U=\U(3,7)$. This map is orientably regular, with orientation-preserving automorphism group $\A^+\M\cong$ $G$, and is fully regular if and only if $K$ is normal in the extended triangle group $\Delta$ of the same type, which contains $\Delta^+$ with index $2$, or equivalently $G$ has an automorphism $\alpha$ (induced by conjugation by $a$, $b$ or $c$) inverting two of the generators $x$, $y$ and $z$.
In this case $\M$ has full automorphism group ${\rm Aut}\,\M\cong\Delta/K$, containing ${\rm Aut}^+\M$ with index $2$, and the elements of ${\rm Aut}\,\M\setminus{\rm Aut}^+\M$ induce anticonformal automorphisms of $\S$.

Every Hurwitz group is perfect, since $\Delta$ is, so it covers a non-abelian finite simple group, which is itself a Hurwitz group. 
Theorem~\ref{th:Macb}, due to Macbeath, gives an important infinite class of simple Hurwitz groups $G={\rm PSL}_2(q)$.

Singerman~\cite{Sin74} has shown that if $q$ is any prime power then each pair of generators of ${\rm PSL}_2(q)$ are simultaneously inverted by some automorphism $\alpha$. It follows that the orientably regular Macbeath--Hurwitz maps $\M$ corresponding to the surfaces $\S$ in Theorem~\ref{th:Macb} are all fully regular. Moreover, since ${\rm Aut}\,{\rm PSL}_2(q)={\rm P\Gamma L}_2(q)$, with $q$ an odd power of $p$ in each case, it follows that the full automorphism group ${\rm Aut}\,\M$ of $\M$ must be either $G\times{\rm C}_2$ or ${\rm PGL}_2(q)$, as the automorphism $\alpha$ of $G$ is respectively inner or outer (this is independent of which pair of generators $x$, $y$ or $z$ it inverts).

More generally, we will call {\em any\/} fully regular orientable map, such as $\M$, {\sl inner} or {\sl outer regular} in these two cases.
If $\M$ is inner regular, it has a non-orientable quotient ${\mathcal N}=\M/{\rm C}_2$ with ${\rm Aut}\,{\mathcal N}\cong G$, and $\M$ is isomorphic to the canonical orientable double cover of $\mathcal N$. Conversely, every non-orientable regular map $\mathcal N$ with automorphism group $G$ has an orientable double cover $\M$ with orientation-preserving automorphism group $G$ and full automorphism group $G\times{\rm C}_2$, so that $\M$ is inner regular. This gives a bijection between inner regular maps $\M$ and non-orientable regular maps~$\mathcal N$.

In this bijection, if $\M$ has genus $g$ and extended type $\{m,n\}_l$ (meaning that it has type $\{m,n\}$ and Petrie length $l$), then ${\mathcal N}$ has genus $g+1$ and extended type $\{m,n\}_{l'}$, where $l'={\rm lcm}\{l,2\}$. Given $\M$, this information about $\mathcal N$ is often sufficient to identify it uniquely in the catalogues, as we will do without further comment in many of our examples. However, there are a few cases where further work is required, and here one can use a MAGMA computation by Conder~\cite{Con25}, which identifies the orientable double covers of all 9982 non-orientable regular maps of genus up to 1502. (Moreover, this list provides a check on all pairings which have been found `by hand' as above.)

In particular, if the map $\M$ of type $\{3,7\}$ and genus $g$ corresponding to a Hurwitz group $G$ is inner regular, then its non-orientable quotient $\N$, having the same type and having genus $g'=g+1$, has automorphism group $G$ attaining  the upper bound $84(g'-2)$ for the number of automorphisms of a non-orientable Klein surface or regular map of genus $g'$, so that $G$ is an $H^*$-group. Conversely, all $H^*$-groups, together with their associated non-orientable maps and Klein surfaces, arise in this way from inner regular Hurwitz groups, acting on their orientable double covers.

\medskip

\noindent{\bf Warning} Not all Hurwitz maps are fully regular: for instance, all of those associated with the `small' Ree groups $^2G_2(3^e)$ for odd $e>1$ are chiral~\cite{Jones:Ree}. The smallest examples are the chiral pair of genus $17$ denoted by C17.1 in~\cite{Con09, Pot}; these can be formed by applying the Macbeath trick, with $m=2$, to Klein's map to give a fully regular 
Hurwitz map of genus $129$, and then factoring out either of the two $3$-dimensional direct summands, conjugate in $\Delta[7,2,3]$, into which the mod~$(2)$ homology module decomposes (see~\cite[Section 6]{JZ} for details).

\medskip

Wendy Hall~\cite[Theorem~2.9]{Hall77} has shown that a Macbeath--Hurwitz map $\M$ is inner regular if and only if $3-t^2$ is a square in $\F_q$, where the canonical generator $y$ of $G={\rm PSL}_2(q)$ of order $7$ corresponding to $\M$ has trace $\pm t$ (see Theorem~\ref{th:Hall}). We will prove a stronger form of this result, namely Theorem~\ref{th:extcriterion}, in Section~\ref{sec:Hall}, but before that we will investigate the consequences of this for the distribution of inner and outer regular Macbeath--Hurwitz maps.

\medskip

\noindent{\bf Example 1} If $q=7$ then $t=\pm 2$, so $3-t^2=-1$ which is not a square in $\F_7$. The corresponding map $\M$ (Klein's map of genus~3, denoted by R3.1 in~\cite{Con09, Pot}) is outer regular, with full automorphism group ${\rm PGL}_2(7)$ and no non-orientable regular quotient. The Riemann surface $\mathcal S$, equivalently the quartic curve
\[x^3y+y^3z+z^3x=0\]
introduced by Klein in~\cite{Kle78}, is examined in detail in the book~\cite{Lev}. (Note that the `obvious' central involution $(x,y,z)\mapsto(-x,-y,-z)$ is in fact the identity, since $x$, $y$ and $z$ are homogeneous coordinates.)

\medskip

\noindent{\bf Example 2} If $q=8$ then since every element of $\F_8$ is a square, the corresponding map R7.1 of genus $7$ has full automorphism group ${\rm PSL}_2(8)\times{\rm C}_2$, with a non-orientable regular quotient N8.1, the only non-orientable regular map of genus $8$ and type $\{3,7\}$. Macbeath described the Riemann surface $\mathcal S$ in~\cite{Macb65}, giving defining equations for the corresponding curve.

\medskip

\noindent{\bf Example 3} If $q=13$ then $G$ has three conjugacy classes of elements of order $7$, with traces $t_j=\pm 6$, $\pm 5$ and $\pm 3$ for $j=1$, $2$ and $3$.
These correspond to the three regular maps R14.1, R14.2 and R14.3 of genus $14$, distinguished by having Petrie lengths $12$, $26$ and $14$.
The corresponding values of $s_j=3-t_j^2$ are $6$, $4$ and $7$, with only $4$ a square in ${\mathbb F}_{13}$, so of these maps, only R14.2 is inner regular, with non-orientable regular quotient N15.1 of Petrie length~$13$.

\medskip

In view of Example~1 we will generally assume, unless otherwise stated, that $q\ne 7$, so that we are in case~(2) or (3) of Theorem~\ref{th:Macb}.


\section{Inner and outer regularity}\label{sec:regularity}

In cases~(2) and (3) of Theorem~\ref{th:Macb} the possible values of $t$ are
\[t_j=t_{-j}:=\zeta^j+\zeta^{-j}\quad (j=1, 2, 3),\]
where $\zeta$ is a primitive $7$th root of unity, contained in $\F_q$ or $\F_{q^2}$ as $q\equiv \pm 1$ mod~$(7)$ respectively. These traces $t_j$ correspond to three Hurwitz maps $\M_j$; these are isomorphic in case~(3), where ${\rm Gal}\,\F_q\cong{\rm C}_3$ induces automorphisms of $G$ permuting the three conjugacy classes of elements of order $7$ transitively, whereas they are mutually non-isomorphic in case~(2), where no such automorphisms exist.

Simple algebra shows that the symmetric functions of the traces $t_j$ are
\[t_1+t_2+t_3=-1,\quad t_1t_2+t_2t_3+t_3t_1=-2\quad\hbox{and}\quad t_1t_2t_3=1,\]
so these traces are the roots in ${\mathbb F}_q$ of the polynomial
\begin{equation}\label{eq:t}
f_0(t):=t^3+t^2-2t-1.
\end{equation}

Let us now define $s_j:=3-t_j^2$, so that $\M_j$ is inner or outer regular as $s_j$ is or is not a square in ${\mathbb F}_q$. Note that
\[s_j=3-t_j^2=3-(\zeta^j+\zeta^{-j})^2=1-(\zeta^{2j}+\zeta^{-2j})=1-t_{2j},\]
with subscripts taken mod~$(7)$. The transformation $j\mapsto 2j$ induces a permutation of the three traces $t_j$. It follows that the symmetric functions of the parameters $s_j$ are
\[s_1+s_2+s_3=4,\quad s_1s_2+s_2s_3+s_3s_1=3\quad\hbox{and}\quad s_1s_2s_3=-1,\]
so these parameters are the roots in ${\mathbb F}_q$ of the polynomial
\begin {equation}\label{eq:s}
f_1(s)=s^3-4s^2+3s+1.
\end{equation}
(Equivalently, we could define $f_1(s)=-f_0(1-s)$, with a minus sign to make this polynomial monic.)

Now let $\chi$ denote the quadratic residue character of ${\mathbb F}_q$, an epimorphism ${\mathbb F}_q^*\to\{\pm 1\}\cong{\rm C}_2$ taking the values $1$ and $-1$ on the squares and non-squares respectively. Since $s_1s_2s_3=-1$ we have
\[\chi(s_1)\chi(s_2)\chi(s_3)=\chi(-1)=\pm 1\quad\hbox{as}\quad q\equiv\pm 1\;{\rm mod}~(4).\]
Thus an odd or even number of the maps $\M_j$ are inner regular as $q\equiv 1$ or $-1$ mod~$(4)$. In case~(3) these three maps are isomorphic, so they have the same regularity properties. We have therefore proved the following:

\begin{proposition}[Theorem 2.10 of \cite{Hall77}]\label{prop:k}
In case~(2), if $p\equiv 1$ {\rm mod}~$(4)$ then either all three maps $\M_j$ are inner regular, or just one of them is, whereas if $p\equiv -1$ {\rm mod}~$(4)$ then either two of the maps $\M_j$ are inner regular, or none of them is.

In case~(3) the unique map $\M$ is inner or outer regular as $p\equiv\pm 1$ {\rm mod}~$(4)$ respectively.
\end{proposition}

The first part of this result is a reformulation of Theorem~\ref{th:parity}.

\medskip

\noindent{\bf Example 4} If $p=43$ we obtain three maps $\M_j$ of genus $474$. (These are R474.1--3 in~\cite{Pot}, with Petrie lengths 22, 42 and 44.) Taking $\zeta=4$ as a primitive $7$th root of $1$ we obtain traces $t_j=15$, $8$ and $19$ for $j=1$, $2$ and $3$, with $s_j=36$, $25$ and $29$. The first two of these are clearly squares, whereas quadratic reciprocity shows that the last is not, so only $\M_1$ and $\M_2$ are inner regular. (Their non-orientable quotients are respectively N475.2 and N475.1, with Petrie lengths 22 and 21.)

\medskip

In case~(3) we thus have a simple criterion for inner or outer regularity; moreover, by de la Vall\'ee-Poussin's quantified form of Dirichlet's Theorem, the sets of primes $p$ giving these two properties have the same density. From now on we will therefore restrict our attention to case~(2), where the situation is less clear.

Examples in~\cite[Ch.~2]{Hall77} with $p=167$, $13$, $43$ and $181$ show that in case~(2) the number $k=k_p$ of inner regular maps $\M_j$ for a given prime $p$ can be $0, 1, 2$ or $3$. For each such $k$ let is define
\[\Sigma_k=\{p\equiv\pm 1\,{\rm mod}~(7)\mid k_p=k\},\]
so that these four sets form a partition of the set of primes $p\equiv\pm 1$ mod~(7).
Conder, extending earlier computations by Hall, has considered the first $400$ such primes, determining how many of them give specific values of $k$ (these primes are listed in the Appendix):

\begin{itemize}
\item for $p\equiv 1$ mod~$(4)$ there are $154$ such primes in $\Sigma_1$, and $47$ in $\Sigma_3$, whereas there are none in $\Sigma_0$ or $\Sigma_2$;
\item for $p\equiv -1$ mod~$(4)$ there are $48$ such primes in $\Sigma_0$, and $151$ in $\Sigma_2$, whereas there are none in $\Sigma_1$ or $\Sigma_3$.
\end{itemize}
Of course, the absence of primes giving even or odd values of $k$ in these two cases is explained by Proposition~\ref{prop:k}. Note that the proportions $47/400$ and $48/400$ are close to $1/8$, while $154/400$ and $151/400$ are close to $3/8$, the asymptotic proportions one would expect if the elements $s_j$ were chosen randomly and independently, with uniform distributions. We will show that these are in fact the relative densities of these sets of primes, as claimed in Theorem~\ref{th:densities}..


\section{Density theorems}\label{sec:density}

This section will explain how the relative densities stated in Theorem~\ref{th:densities} can be obtained from the density theorems of Frobenius and of Chebotarev (see~\cite{LS} for a very readable account of the life and work of the latter, and the relationship between these two theorems).

A root $s_j$ of $f_1(s)$ is a square $s_j=a_j^2$ in ${\mathbb F}_p$ if and only if $a_j$ (together with $-a_j$) is a root of the polynomial
\begin{equation}
f_2(x):=f_1(x^2)=x^6-4x^4+3x^2+1
\end{equation}
in ${\mathbb F}_p$, or equivalently, $f_2(x)$ has a pair of linear factors $x\pm a_j$ in ${\mathbb F}_p[x]$. The number $k_p$ of inner regular maps $\M_j$ for each $p$ is therefore twice the number of linear factors of $f_2(x)$ in ${\mathbb F}_p[x]$.

\medskip

\noindent{\bf Example 5} If $p=13$ the factorisation is
\[f_2(x)=(x-2)(x+2)(x^2-6)(x^2-7),\]
with a pair of linear factors $x\pm 2$ corresponding to the root $s_j=4=(\pm 2)^2$ of $f_1$ and thus to the unique inner regular map $\M_j$, and two irreducible quadratic factors corresponding to the non-square roots $6$ and $7$ and thus to the two outer regular maps.

\medskip

Any polynomial $f(x)\in{\mathbb Z}[x]$, when reduced mod~$(p)$ for some prime $p$, has a factorisation as a product of irreducible polynomials in $\F_p[x]$. If $f$ is monic then the degrees of these factors form a partition $\pi$ of its degree $d=\deg(f)$, which depends on $p$. If the discriminant of $f$ is non-zero, then $f$ has $d$ simple roots, permuted faithfully by its Galois group as a subgroup of the symmetric group ${\rm S}_d$. Frobenius's theorem~\cite{Fro} (proved, according to~\cite{LS}, in 1880, though not published until 1896) states that the density of the set of primes giving a particular partition $\pi$ is equal to the proportion of elements of the Galois group which have cycle-lengths (as permutations of the roots) giving the partition $\pi$.

Before proceeding, we need to explain more precisely what we mean here by density.
Let $\Pi$ denote the set of all prime numbers, and let us define
\[\Pi'=\{p\in\Pi\mid p\equiv \pm 1\, {\rm mod}~(7)\},\]
\[\Pi''=\{p\in\Pi\mid p\equiv \pm 2 \;\hbox{or}\;\pm 3\, {\rm mod}~(7)\}.\]
Thus $\Pi$ has a partition
\[\Pi=\{7\}\,\dot{\cup}\,\Pi'\,\dot{\cup}\,\Pi''\]
corresponding to the three cases in Theorem~\ref{th:Macb}.
We will say that a subset $\Sigma\subseteq\Pi$ has {\sl natural density\/} $d$ (within $\Pi$) if
\[\lim_{x\to +\infty}\frac{|\Sigma\cap [0,x]|}{|\Pi\cap [0,x]|}=d.\]
For example, it follows from the quantified form of Dirichlet's theorem on primes in arithmetic progressions that $\Pi'$ and $\Pi''$ have natural densities $2/\varphi(7)=1/3$ and $4/\varphi(7)=2/3$.
We will say that a subset $\Sigma\subseteq\Pi'$ has {\sl relative density\/} $d'$ (within $\Pi'$) if
\[\lim_{x\to +\infty}\frac{|\Sigma\cap [0,x]|}{|\Pi'\cap [0,x]|}=d'.\]
Then $d'=3d$ if these densities exist.

\medskip

\noindent{\bf Example 6} An irreducible cubic polynomial in ${\mathbb Z}[x]$ has Galois group ${\rm A}_3$ or ${\rm S}_3$ as its discriminant is or is not a square. A general cubic $ax^3+bx^s+cx+d$ has discriminant $b^2c^2-4ac^3-4b^3d-27a^2d^2+18abcd$, so the polynomial $f_1(x)=x^3-4x^2+3x+1$ in (\ref{eq:t}) has discriminant $49=7^2$; it therefore has Galois group $G={\rm A}_3\cong{\rm C}_3$, permuting the three complex (in fact real) roots regularly. One third of the elements of $G$ (namely the identity) have cycle structure $1^3$, so the primes for which $f_1$ splits into linear factors (i.e.~has three roots in ${\mathbb F}_p$) have density $1/3$; these are, of course, the primes $p\equiv\pm 1$ mod~$(7)$ forming $\Pi'$, with roots $s_j\in{\mathbb F}_p$. The other two elements of $G$ induce $3$-cycles, so the primes for which $f_1$ is irreducible have density $2/3$; these are the primes $p\equiv \pm 2$ or $\pm 3$ mod~$(7)$ forming $\Pi''$. (If $p=7$ then $f_1$ reduces to $(x+1)^3$, with a triple root at $-1$; this `bad reduction' is reflected in the discriminant being divisible by $7$.)

\medskip

Chebotarev's theorem~\cite{Che23,Che25} refines that of Frobenius by associating sets of primes, not with cycle-structures but with specific conjugacy classes in the Galois group, a finer distinction. (Of course, conjugate elements have the same cycle-structure, but the converse may fail.) It also considers  arbitrary Galois extensions of algebraic number fields, not just extensions of $\mathbb Q$. Without going into the general theory (see~\cite{LS} for this), it is sufficient here to explain how it works in our particular situation.


\section{Application to $f_1$}\label{sec:f1}

The splitting field of $f_1$ is a Galois extension of $\mathbb Q$ of degree $|G_1|=3$. This is the field
\[K_1:={\mathbb Q}(\tau_j)={\mathbb Q}(\sigma_j)={\mathbb Q}(\zeta_7)\cap{\mathbb R}\]
obtained by adjoining any complex root $\tau_j$ of $f_0$ or $\sigma_j$ of $f_1$ to $\mathbb Q$, equal to
the real subfield of the cyclotomic field ${\mathbb Q}(\zeta_7)$ where $\zeta_7=\exp(2\pi i/7)$. The Galois group $G_1$ of $f_1$ is generated by the automorphism $\theta:\sigma_1\mapsto\sigma_2$ of order $3$ of $K_1$ induced by squaring each complex $7$th root of $1$. This permutes the roots $\tau_j$ of $f_0$, and hence the roots $\sigma_j=3-\tau_j^2$ of $f_1$, by inducing the $3$-cycle $(1,2,3)$ on their subscripts.

The Frobenius automorphism $\phi_p$ acts on the algebraic closure $\overline{\mathbb F}_p$ of ${\mathbb F}_p$ by raising every element to its $p$th power. It fixes ${\mathbb F}_p$, and hence induces a permutation of the roots of any polynomial over ${\mathbb F}_p$. In the case of $f_1$, if $p\equiv\pm 2$ mod~$(7)$ then since $\phi_p$ acts by squaring (and possibly inverting) the $7$th roots of $1$ in $\overline{\mathbb F}_p$, it also induces the permutation $(1,2,3)$ of the subscripts of the roots $s_j$ of $f_1$ in this field, whereas if $p\equiv\pm 3$ mod~$(7)$ it induces the inverse permutation. Thus the two sets of primes $p\equiv\pm 2$ and $\pm 3$ mod~$(7)$ forming $\Pi''$ are associated with these two permutations in $G_1$ (more correctly, with their conjugacy classes $C$ in $G_1$, which are singletons in this abelian group), and their natural densities are equal to the relative sizes $|C|/|G_1|$ in $G_1$ of these classes, namely $1/3$ in each case. Chebotarev's theorem tells us more here than that of Frobenius: these two permutations have the same cycle-structure but are not conjugate, and they are each associated with an explicit set of primes. For completeness, we note that if $p\equiv\pm 1$ mod~$(7)$ then $\phi_p$ induces the identity permutation on the roots $s_j$, and the identity class in $G_1$ also has relative size $1/3$, equal to the natural density of $\Pi'$.


\section{Application to $f_2$}\label{sec:density}

To solve our problem we need to understand the factorisation of the reduction mod~$(p)$ of $f_2$, rather than $f_1$. If we are to apply density theorems to $f_2$ we need to know its Galois group $G_2$.
Since $f_2$ is the composition of two polynomials, namely $x\mapsto x^2$ and $f_1$, this is an imprimitive group of degree $d=\deg(f_2)=6$, permuting the three pairs of roots $\pm \alpha_j$ ($j=1,2,3$) in $\mathbb C$ in the same way as the Galois group $G_1={\rm A}_3\cong{\rm C}_3$ of $f_1$ permutes its roots $\sigma_j=\alpha_j^2$. Thus $G_2$ is a subgroup of the wreath product
\[W:={\rm S}_2\wr{\rm A}_3\cong{\rm C}_2\wr{\rm C}_3\cong{\rm V}_8\rtimes{\rm C}_3,\]
a semidirect product of an elementary abelian normal subgroup $N\cong{\rm V}_8={\rm C}_2\times{\rm C}_2\times{\rm C}_2$ by a complement $H\cong{\rm C}_3$ which transitively permutes the direct factors ${\rm C}_2$ by conjugation. By considering sign-changes of $f_1$ over the interval $[-1,3]\subset{\mathbb R}$, we see that two of its three roots $\sigma_j$ are positive and one (namely $\sigma_3$) is negative, so $f_2$ has four real roots and a complex conjugate imaginary pair $\pm\alpha_3$. Thus $G_2$ contains an involution, namely complex conjugation, transposing this last pair and fixing the other four roots; the conjugates of this involution under $H$ generate the subgroup $N$ of $W$, so $G_2=W$.

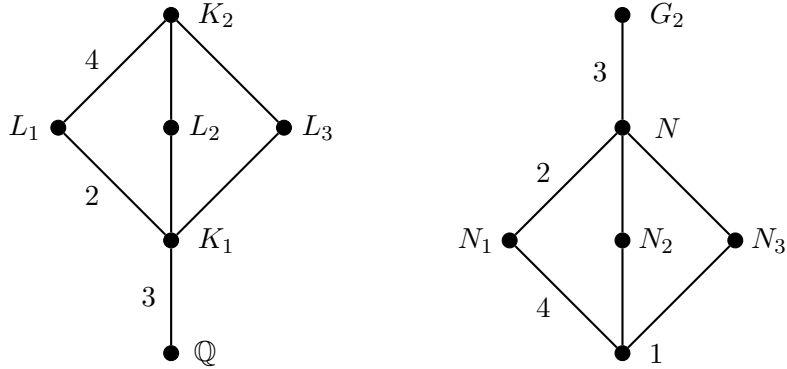
\begin{figure}[h!]
\begin{center}

 \begin{tikzpicture}[scale=0.3, inner sep=0.7mm]
 
\node (A) at (-10,0)  [shape=circle, draw, fill=black] {};
\node (B) at (-10,5)  [shape=circle, draw, fill=black] {};
\node (C) at (-15,10)  [shape=circle, draw, fill=black] {};
\node (D) at (-10,10)  [shape=circle, draw, fill=black] {};
\node (E) at (-5,10)  [shape=circle, draw, fill=black] {};
\node (F) at (-10,15)  [shape=circle, draw, fill=black] {};

\draw [thick] (A) to (B) to (D) to (F);
\draw [thick] (F) to (C) to (B) to (E) to (F);

\node at (-8,15) {$K_2$};
\node at (-8,5) {$K_1$};
\node at (-8.5,0) {$\mathbb Q$};
\node at (-16.5,10) {$L_1$};
\node at (-8.5,10) {$L_2$};
\node at (-3.5,10) {$L_3$};

\node at (-13.5,13) {$4$};
\node at (-13.5,7) {$2$};
\node at (-11,2.5) {$3$};


\node (a) at (10,15)  [shape=circle, draw, fill=black] {};
\node (b) at (10,10)  [shape=circle, draw, fill=black] {};
\node (c) at (15,5)  [shape=circle, draw, fill=black] {};
\node (d) at (10,5)  [shape=circle, draw, fill=black] {};
\node (e) at (5,5)  [shape=circle, draw, fill=black] {};
\node (f) at (10,0)  [shape=circle, draw, fill=black] {};

\draw [thick] (a) to (b) to (d) to (f);
\draw [thick] (f) to (c) to (b) to (e) to (f);

\node at (11.5,0) {$1$};
\node at (12,10) {$N$};
\node at (12,15) {$G_2$};
\node at (16.5,5) {$N_3$};
\node at (11.5,5) {$N_2$};
\node at (3.5,5) {$N_1$};

\node at (6.5,2) {$4$};
\node at (6.5,8) {$2$};
\node at (9,12.5) {$3$};

\end{tikzpicture}

\end{center}
\caption{Some subfields of $K_2$ and the corresponding subgroups of $G_2$.}
\label{fig:subfields}

\end{figure}

It follows from the Galois correspondence that the splitting field $K_2$ of $f_2$ has degree $|G_2|=24$ over $\mathbb Q$, and has degree $8$ over the splitting field $K_1$ of $f_1$. Adjoining a square root $\alpha_j$ of $\sigma_j$ gives a quadratic extension $L_j$ of $K_1$, and these three subfields $L_j$ generate $K_2$. See Figure~\ref{fig:subfields} for these subfields of $K_2$ and their corresponding subgroups of $G_2$, with $N_j\cong{\rm V}_4$ corresponding to $L_j$; integers represent the degrees and indices of the corresponding inclusions.

We can now list the elements of $G_2$, together with their cycle-structures on the six roots $\pm \alpha_j$ of $f_2$. They are as follows:
\begin{itemize}
\item the identity element, with cycle-structure $1^6$;
\item seven elements of order $2$, three each with cycle-structures $1^42^1$ and $1^22^2$, and one, generating the centre of $G_2$, with cycle-structure $2^3$;
\item eight elements of order $3$, all with cycle-structure $3^2$;
\item eight elements of order $6$, all with cycle-structure $6^1$.
\end{itemize}
It follows from Frobenius's theorem that there are
\begin{itemize}
\item[(a)] two sets of primes $p$, each with density $1/24$, such that the reduction of $f_2$ mod~$(p)$ is respectively a product of six linear factors or three irreducible quadratics;
\item[(b)] two sets, each of density $1/8$, for which it has four or two linear factors and one or two irreducible quadratics;
\item[(c)] two sets, each of density $1/3$, for which it is irreducible or a product of two irreducible cubics.
\end{itemize}
This implies that the sets of primes $p$ for which the reduction of $f_2$ mod~$(p)$ has $k_p=0$, $1$, $2$ or $3$ pairs $x\pm a_j$ of linear factors have natural densities $17/24$, $1/8$, $1/8$ and $1/24$.

In each of cases (a) and (b), the two sets described are conjugacy classes in $G_2$; in (a) they are singletons, the two elements forming the centre of $G_2$, while the sets in (b) form the rest of $N $. However, in (c) each of the two sets described is the union of two conjugacy classes, distinguished by inclusion in different cosets of $N$; this comes from lifting back to $G_2$ the partition of $G_1=G_2/N\cong{\rm C}_3$ into (singleton) conjugacy classes, used in Section~\ref{sec:f1}.

These densities do not match the computational evidence given at the end of Section~\ref{sec:regularity} because they are densities within the set $\Pi$ of {\sl all} primes, whereas the only primes relevant in case~(2) are those $p\equiv\pm 1$ mod~$(7)$ forming $\Pi'$. By Chebotarev's theorem, the primes $p$ in (c), namely those for which there are irreducible factors of degree divisible by $3$, are those for which $\phi_p$ induces a permutation of the roots of order divisible by $3$, corresponding to the cubic fields in case~(3) where $p\in\Pi''$. If we exclude those primes then the remaining set $\Pi'$, with density $1/3$ in $\Pi$, is partitioned into subsets $\Sigma_k$ with {\sl relative\/} densities $1/8$, $3/8$, $3/8$ and $1/8$ in $\Pi'$ for $k=0, 1, 2$ and $3$, as claimed. This proves Theorem~\ref{th:densities}.

Note also that, by Proposition~\ref{prop:k}, $\Sigma_1\cup\Sigma_3$ and $\Sigma_0\cup\Sigma_2$ are the subsets of $\Pi'$ consisting of those primes $p\equiv\pm 1$ mod~$(4)$ respectively. 

\medskip

\noindent{\bf Problem} Can one give necessary {\sl and sufficient} conditions for a prime $p\equiv \pm 1$ mod~$(7)$ to be a member of a particular set $\Sigma_k$?


\section{An alternative approach}\label{sec:alternative}

An alternative way of restricting attention to primes $p\in\Pi'$ is to consider the Galois group of $f_2$, not over $\mathbb Q$, but over $K_1$, where it has the factorisation
\begin{equation}
f_2(x)=(x^2-\sigma_1)(x^2-\sigma_2)(x^2-\sigma_3).
\end{equation}
The motivation for this is that since the fields ${\mathbb F}_p$ appearing in case~(2), with $p\equiv\pm 1$ mod~$(7)$, all contain the traces $t_j=\zeta^j+\zeta^{-j}$ of the elements of order $7$ in ${\rm PSL}_2(p)$, it is appropriate to regard $f_2$ as a polynomial, not over $\mathbb Q$ but over the field $K_1$ generated by the traces $\tau_j=\zeta_7^j+\zeta_7^{-j}$ in $\mathbb C$, or equivalently by the roots $\sigma_j$ of $f_1$. 

Now the triangle group $\Delta^+$ is an arithmetic group, defined over $K_1$. 
In the ring ${\mathcal O}={\mathbb Z} [2\cos 2\pi/7]$ of integers of this field, rational primes $p\in\Pi''$ generate prime ideals $\mathfrak p$ (that is, they are `inert'), with ${\mathcal O}/{\mathfrak p}\cong\F_{p^3}$, whereas any prime $p\in\Pi'$ generates a product of three prime ideals ${\mathfrak p}_j$ ($j=1, 2, 3$), each with ${\mathcal O}/{\mathfrak p}_j\cong\F_p$. For example $13=(1-2\sigma_1)(1-2\sigma_2)(1-2\sigma_3)$. (The prime $7$ ramifies in $\mathcal O$, so its quotient has nilpotent elements and is not a field; for instance $(\sigma_j+1)^3=7\sigma_j^2$, so $(s_j+1)^3=0$ when $p=7$.) As noted briefly by Magnus~\cite[Section II.7]{Mag}, the Macbeath--Hurwitz groups in cases~(2) and (3) are obtained from $\Delta^+$ by factoring out its principle congruence subgroups corresponding to these prime ideals ${\mathfrak p}_j$ and $\mathfrak p$; see~\cite{Dza} for D\v zambi\'c's detailed explanation of this viewpoint.

As a polynomial in $K_1[x]$ with splitting field $K_2$, $f_2$ has Galois group  ${\rm Gal}\,K_2/K_1\cong{\rm V}_8$ consisting of the elements of order $1$ or $2$ in $G_2$ described in Section~\ref{sec:density}. The Frobenius automorphism $\phi_p$ acts on the roots of $f_2$ by
\[\phi_p:a_j\mapsto a_j^p=s_j^{(p-1)/2}a_j=\pm a_j\]
as $s_j^{(p-1)/2}=\pm 1$, that is, as $s_j$ is or is not a square in ${\mathbb F}_p$. It thus induces a permutation of order $1$ or $2$, with pairs of fixed points corresponding to inner regular maps $\M_j$, and transpositions to outer regular maps.

In the group $N$ the proportion of elements with $k=0, 1, 2$ or $3$ pairs of fixed points is $1/8$, $3/8$, $3/8$ and $1/8$ respectively, so by Chebotarev's extension of Frobenius's theorem to arbitrary Galois extensions of algebraic number fields, these are the relative densities of the sets $\Sigma_k$, as claimed.

\medskip

\noindent{\bf Problem} Can one use this approach to be more specific about which primes are in each set $\Sigma_k$?


\section{Macbeath maps of type $\{3,n\}$}\label{sec:Macbeathmaps}

In the rest of this paper we will consider how much of this work extends to the groups, surfaces and maps where the latter have type $\{3,n\}$ for some $n\ge 7$. We will therefore redefine 
\[\Delta^+:=\Delta(2,3,n)=\langle X, Y, Z\mid X^3=Y^n=Z^2=XYZ=1\rangle\]
from now on. The genus $g$ of such an orientably regular map $\M$, with $G=\A^+\M$, is given by
\[|G|=\frac{12n}{n-6}(g-1),\]
so when $G={\rm PSL}_2(q)$, as we will continue to assume, we have
\[g=1+\frac{n-6}{12n}\cdot\frac{q(q^2-1)}{\gcd(q,2)}\]
\[=1+\frac{(n-6)\cdot q(q^2-1)}{24n}
\quad\hbox{or}\quad 1+\frac{(n-6)\cdot q(q^2-1)}{12n}\]
as $q$ is odd or even.

Let us assume that $p$ is coprime to $n$. Then $G$ is a smooth quotient $\Delta^+/K$ of $\Delta^+$ if and only if $q=p^d$ where $d$ is the multiplicative order of $p$ in ${\mathbb Z}_n^*/\{\pm 1\}$. There are $\varphi(n)/2$ conjugacy classes of elements $y$ of order $n$ in $G$, each of them giving, via a generating triple $x$, $y$, $z$, an orientably regular map $\M$ of type $\{3,n\}$, which we call a Macbeath map. Two such maps are isomorphic if and only if the traces of their canonical generators $y$ are conjugate under ${\rm Gal}\,\F_q\cong{\rm C}_d$, so we obtain $\varphi(n)/2d$ isomorphism classes of maps, and also of Riemann surfaces $\S$. (We have ${\rm Aut}\,\S=G$ since $\Delta^+=\Delta(2,3,n)$ is finitely maximal, see~\cite{Sin72}, so that it is the normaliser of $K$ in $\PSL_2({\mathbb R})$.) By Singerman's result on inverting generating pairs for $G$, these maps are all fully regular, with full automorphism group $G\times{\rm C}_2$ or ${\rm PGL}_2(q)$ as the automorphism of $G$ inverting two of the canonical generators $x$, $y$, $z$ is inner or outer. In the next section we will prove Theorem~\ref{th:extcriterion}, that these two cases, inner and outer regularity, arise as $3-t^2$ is or is not a square in $\F_q$, where $y$ has trace $\pm t$.


\section{Extension of Hall's Criterion}\label{sec:Hall}

Here we summarise Hall's proof of Theorem~\ref{th:Hall} in~\cite{Hall77}, omitting any statements (explicit or implicit) that $n=7$ since these are never essential to the argument. This means that the proof in fact applies to Macbeath maps of type $\{3,n\}$ for all $n\ge 7$, as claimed in Theorem~\ref{th:extcriterion}. In some places we have changed Hall's notation to make it consistent with that in the rest of this paper.

\medskip

Identifying each matrix in ${\rm SL}_2(q)$ with its negative, as usual, let
\[z=\begin{pmatrix} 0 & 1 \\ -1  & 0 \end{pmatrix}
\quad\hbox{and}\quad
x=\begin{pmatrix} r & s \\ u & v \end{pmatrix}\]
 be generators of $G={\rm PSL}_2(q)$ of orders $2$ and $3$, with $y=(zx)^{-1}$ of order $n$. (Since all involutions are conjugate in $G$, we may assume that $z$ has this form.) By Singerman's result (Hall proves this by a direct construction) there is a (unique) involution $w\in{\rm PGL}_2(q)$ which, acting by conjugation, inverts $z$ and $x$, that is, $w$, $wz$ and $wx$ are involutions, or equivalently they have trace $0$. If we write $w=\begin{pmatrix} \alpha & \beta \\ \gamma & \delta \end{pmatrix}$ the first and second trace conditions give $\alpha+\delta=-\beta+\gamma=0$, so that $w=\begin{pmatrix} \alpha & \beta \\\beta & -\alpha \end{pmatrix}$. The corresponding map $\M$ is inner regular if and only if $w\in{\rm PSL}_2(q)$, or equivalently $\det w=-\alpha^2-\beta^2$ is a (necessarily non-zero) square in ${\mathbb F}_q$.

The third trace condition gives
$\alpha r+\beta u+\beta s-\alpha v=0$, so that
\begin{equation}\label{eq:alphabeta}
\alpha^2(r-v)^2=\beta^2(s+u)^2.
\end{equation}
(Squaring here introduces a second involution $w'$, which inverts $z$ but not $y$; however, everything which follows applies to $w$.)
Adding $\beta^2(r-v)^2$ to each side in (\ref{eq:alphabeta}), and then expanding the squares on the right-hand side gives
\begin{equation}\label{eq:squares}
(\alpha^2+\beta^2)(r-v)^2=\beta^2(r^2+s^2+u^2+v^2-2)
\end{equation}
since $rv-su=\det x=1$.

Macbeath~\cite{Macb69} showed that if any pair $z$ and $x$ generate ${\rm PSL}_2(q)$, then the quadratic form
\[Q(\xi, \eta, \zeta)
=\xi^2+\eta^2+\zeta^2+a\eta\zeta+b\zeta\xi+c\xi\eta,\]
where $a, b, c$ are the traces of $z$, $x$ and $zx$,
is nonsingular, which implies that the determinant
\begin{equation}\label{eq:det}
\det
\begin{pmatrix}
1 & c/2 & b/2 \\
c/2 & 1 & a/2 \\
b/2 & a/2 & 1 \\
\end{pmatrix}
=
\frac{1}{4}(4+abc-a^2-b^2-c^2)
\end{equation}
is non-zero. In our case we have $a=0$, $b=r+v$ and $c=u-s$, and from this it follows that $r^2+s^2+u^2+v^2\ne 2$. Thus the second factor on the right-hand side of (\ref{eq:squares}) is non-zero.

Now $(r+v)^2={\rm tr}^2\, x=1$ and $(u-s)^2={\rm tr}^2\,zx=t^2$, so adding and using $\det x=1$ again gives
\[1+t^2=r^2+v^2+2rv+u^2+s^2-2su=r^2+s^2+u^2+v^2+2.\] 
and hence
\[t^2-3=r^2+s^2+u^2+v^2-2,\]
with both sides non-zero, as shown above. Substituting this in the right-hand side of (\ref{eq:squares}) gives
\begin{equation}\label{eq:squares2}
(\alpha^2+\beta^2)(r-v)^2=\beta^2(t^2-3).\end{equation}
Provided both sides of (\ref{eq:squares2}) are non-zero (a possibility not considered in~\cite{Hall77}), we see that $-(\alpha^2+\beta^2)$ is a square if and only if $3-t^2$ is, as required.

What if both sides of (\ref{eq:squares2}) are zero? Clearly $\alpha^2+\beta^2=-\det w\ne 0$, and we saw above that $t^2-3\ne 0$, so the only possibility is that $\beta=0$ and $r=v$. Then $\det w=-\alpha^2$ is a square if and only if $-1$ is. If $r=v$ then since ${\rm tr}\, x=1$ we have $r=v=1/2$ and hence $su=rv-\det y=-3/4$, while $-s+u={\rm tr}\,zx=t$, so $-s$ and $u$ are the roots of the quadratic polynomial $T^2-tT+3/4$, the discriminant $t^2-3$ of which must be a square. Thus $3-t^2$ is a square if and only if $-1$ is, the same condition as for $\det w$ to be a square.

 This completes the proof of Theorem~\ref{th:Hall}, and also, since the condition that $n=7$ has not been used, that of Theorem~\ref{th:extcriterion}.
\hfill$\square$


\section{Parity}\label{sec:parity}

In the case $n=7$, Hall~\cite{Hall77} used the quadratic residue character of ${\mathbb F}_p$ to determine the parities of the numbers $k_p$ and $l_p$ of inner and outer regular Macbeath--Hurwitz maps associated with a particular prime $p$, as stated in Theorem~\ref{th:parity}. Here we will show how to extend this to other values of $n$.

Since $k_p+l_p=\varphi(n)/2$, for any given $n$ it is sufficient to determine the parity of $l_p$. We will show that it is possible to evaluate $\prod_js_j$, and hence to find $\prod_j\chi(s_j)=\chi(\prod_js_j)$, where $\chi$ is the quadratic residue character of ${\mathbb F}_q$. Since $\M_j$ is inner or outer regular as $\chi(s_j)=\pm 1$, this tells us the parity of the number of $j$ with $\M_j$ outer regular. Now the $\varphi(n)/2$ maps $\M_j$ are grouped, by the action of ${\rm Gal}\,{\mathbb F}_q\cong{\rm C}_d$, into $\varphi(n)/2d$ sets of $d$ mutually isomorphic maps, so if $d$ is even we obtain no new information about $k_p$ and $l_p$ since they count isomorphism classes of maps: the product is always $1$, with factors $\chi(s_j)=-1$ occurring in even multiples. However, if $d$ is odd we can ignore repetitions, and obtain the parity of the number $l_p$ of non-isomorphic outer regular maps $\M_j$, and hence that of $k_p$.

We will concentrate mainly on the primes $p\equiv\pm 1$ mod~$(n)$, for which $d=1$, so that there are $\varphi(n)/2$ maps with ${\rm Aut}^+\M_j\cong{\rm PSL}_2(p)$. For simplicity we will first deal with the case where $n$ is odd. .The analogue of the cubic polynomial $f_0(t)$ with roots $t_j$ is the minimal polynomial $\Psi_n(t)$ of $\zeta_n+\zeta_n^{-1}=2\cos(2\pi/n)$. (For instance, $\Psi_7(t)=t^3+t^2-2t-1=f_0(t)$.) As before $s_j=1-t_{2j}$, with subscripts mod~$(n)$, so if $n$ is odd the elements $s_j$ are the roots of the monic polynomial \[f_1(s)=(-1)^{\varphi(n)/2}\Psi_n(1-s),\]
which is irreducible since $\Psi_n(t)$ is. In the case $n=7$ the quadratic residue character $\chi$ was applied to $\prod_j s_j$ in order to determine the parities of the numbers $k_p$ and $l_p$ of inner and outer regular maps for different primes $p$. This can be done for more general values of $n$, as follows.

Since $f_1$ has degree $\varphi(n)/2$ we have
\[\prod_js_j=(-1)^{\deg f_1}f_1(0)=\Psi_n(1),\]
so $\prod_j\chi(s_j)=\chi(\Psi_n(1))$. Now Watkins and Zeitlin~\cite{WZ} have shown that for odd $n=2m+1$ we have
\[\prod_{e\mid n}\Psi_e(2x)=2(T_{m+1}(x)-T_m(x)),\]
where $T_m$ is the $m$th Chebyshev polynomial of the first kind. Putting $x=\frac{1}{2}=\cos(\pi/3)$ we have
\[2(T_{m+1}(x)-T_m(x)=c_m:=2\bigl(\cos\bigl(\frac{(m+1)\pi}{3}\bigr)-\cos\bigl(\frac{m\pi}{3}\bigr)\bigr),\]
by the trigonometric definition of the Chebyshev polynomials $T_m$. Now $\cos(m\pi/3)=1$, $1/2$, $-1/2$ or $-1$ as $m\equiv 0$, $\pm 1$, $\pm 2$ or $3$ mod~$(6)$, so
\[c_m=-1,\, -2,\, -1,\, 1,\, 2\;\;\hbox{or}\;\; 1\]
as $m\equiv 0, 1, 2, 3, 4$ or $5$ mod~$(6)$, and hence
\[\prod_{e\mid n}\Psi_e(1)=b_n:=c_m=-1,\, -2,\, -1,\, 1,\, 2\;\;\hbox{or}\;\; 1\]
as $n=2m+1\equiv 1,\, 3,\, 5,\, 7,\, 9$ or $11$ mod~$(12)$.
M\"obius inversion therefore gives
\[\Psi_n(1)=\prod_{e\mid n}b_e^{\,\mu(n/e)}\]
where $\mu$ is the classical M\"obius function, so that
\[\prod_j\chi(s_j)=\chi(\Psi_n(1))=\prod_{e\mid n}\chi(b_e)^{\mu(n/e)}.\]

In order to evaluate the right-hand side we need to know the values of $\chi$ at $b_e=
\pm 1, \pm 2$.
\begin{itemize}
\item If $e\equiv 7$ or $11$ mod~$(12)$ then $\chi(b_e)=\chi(1)=1$.
\item If $e\equiv 1$ or $5$ mod~$(12)$ then $\chi(b_e)=\chi(-1)=\pm 1$ as $p\equiv\pm 1$ mod~$(4)$.
\item If $e\equiv 9$ mod~$(12)$ then $\chi(b_e)=\chi(2)=\pm 1$ as $p\equiv 1, -1
$ mod~$(8)$ or $p\equiv 3,-3$ mod~$(8)$.
\item If  $e\equiv 3$ mod~$(12)$ then $\chi(b_e)=\chi(-2)=\pm 1$ as $p\equiv 1, 3
$ mod~$(8)$ or $p\equiv -1,-3$ mod~$(8)$.
\end{itemize}
Knowing the right-hand side, one can determine the parity of the number $l_p$ of outer regular maps $\M_j$ (those with $\chi(s_j)=-1$) for each prime $p$, and hence that of the number $k_p$ of inner regular maps, thus giving an analogue of Theorem~\ref{th:parity} for each $n$. Some illustrative examples of this will be given later.

When $n$ is even we use traces $t_j=\zeta_{2n}^j+\zeta_{2n}^{-j}$ instead, but again we have $s_j=3-t_j^2=1-t_{2j}$, so that $\prod_js_j=\Psi_n(1)$, as when $n$ is odd. For even $n=2m$  we have
\[\prod_{e\mid n}\Psi_e(2x)=2(T_{m+1}(x)-T_{m-1}(x)).\]
Again taking $x=\frac{1}{2}=\cos(\pi/3)$, we now have
\[2(T_{m+1}(x)-T_{m-1}(x))=c'_m:=2\bigl(\cos\bigl(\frac{(m+1)\pi}{3}\bigr)-\cos\bigl(\frac{(m-1)\pi}{3}\bigr)\bigr),\]
which takes the values
\[c'_m=0,\, -3/2,\, -3/2,\, 0,\, 3/2,\;\hbox{or}\; 3/2\]
as $m\equiv 0, 1, 2, 3, 4$ or $5$ mod~$(6)$. Thus for even $n$ we have
\[\prod_{e\mid n}\Psi_e(1)=b_n:=c_m=0,\, -3/2,\, -3/2,\, 0,\, 3/2\;\hbox{or}\;\; 3/2\]
as $n=2m\equiv 0,\, 2,\, 4,\, 6,\, 8$ or $10$ mod~$(12)$. It again follows by M\"obius inversion that
\[\Psi_n(1)=\prod_{e\mid n}b_e^{\,\mu(n/e)}\]
but with these new values of $b_e$ for even $e$. For general interest, and for use in later examples, the values of $\Psi_n(1)$ for some small $n\ge 7$ are presented in the following table.

\begin{table}[htbp]
\begin{tabular}{c||c|c|c|c|c|c|c|c|c|c|c|c|c}
$n$ &7 & 8 & 9 & 10 &11 & 12 & 13 & 14 & 15 & 16 & 17 & 18 & 19 \\
\hline
$\Psi_n(1)$ & $-1$ &  $-1$ & $-1$ & $-1$ & $-1$ & $-2$ & $1$ & $-1$ 
 & $1$ &  $-1$ &  $1$ & $-3$ & $-1$  \\

\end{tabular}
\vspace{2mm}
\caption{The values of $\Psi_n(1)$ for $n=7,\ldots, 19$.}
\label{tab:Psivalues}
\end{table}


\section{Examples with $n$ odd}\label{sec:oddexamples}

Here we give illustrative examples for some small odd $n>7$. (Even values are considered in the following section, after an explanation of some small technical complications involved in them.) We will deal with the cases $n=9$ and $11$ in some detail, and then give briefer summaries for the odd integers $n=13,\ldots, 19$. For each such $n$ we will consider a few primes $p$, usually chosen small enough that we can compare our results with the evidence in the catalogues compiled by Conder and Poto\v cnik~\cite{Con09, Pot}, the latter listing the regular maps for genus $g=2,\ldots,1501$.

For each pair $n$ and $p$ one can compute the multiplicative order $d$ of $p$, giving the group $G={\rm PSL}_2(q)$ where $q=p^d$, and hence the genus $g$ of the corresponding maps. The number $k_p+l_p$ of them can be checked by counting entries R$g.i$ of type $\{3,n\}$ in the catalogues. Each such map is inner regular if and only if it has a non-orientable quotient of the same type and of genus $g+1$, so $k_p$ and hence $l_p$ can be found from the catalogues. In all cases, one can check that $l_p$ satisfies the appropriate parity conditions, given by evaluating $\Psi_n(1)$. In most cases, different maps R$g.i$, for a given type $\{3,n\}$ and genus $g$, are distinguished from each other by their Petrie lengths $r$, which are given in the catalogues as part of the extended type $\{3,n\}_r$ of these maps; any non-orientable quotient must have Petrie length $r$ equal to or half of that of its double cover, as $r$ is even or odd, so this often allows one to identify which of the maps R$g.i$ are inner regular, and which are not.

In many cases, where $q$ is sufficiently small, Hall's criterion has been checked by computing the trace $\pm t$ corresponding to each map, and determining whether or not $3-t^2$ is a square in ${\mathbb F}_q$, for instance by using quadratic reciprocity in cases where $q=p$. This gives an alternative method of finding $k_p$ and $l_p$; in all cases where this has been done, the results match those obtained as above from the catalogues.

\medskip

\noindent{\bf Example 7} Let $n=9$, so that $\varphi(n)/2=3$, as in the case $n=7$ considered earlier. We have a similar situation, that if $p\equiv\pm 1$ mod~$(9)$ there are three non-isomorphic maps, with $d=1$, whereas if $p\equiv\pm 2$ or $\pm 4$ mod~$(9)$ there are three mutually isomorphic maps, with $d=3$. We have
\[\Psi_1(2x)\Psi_3(2x)\Psi_9(2x)=2(T_5(x)-T_4(x))\]
\[=2((16x^5-20x^3+5x)-(8x^4-8x^2+1))\]
\[=(2x-2)(2x+1)(8x^3-6x+1),\]
so that
\[\Psi_9(x)=x^3-3x+1,\]
and thus $\Psi_9(1)=-1$. Alternatively,
\[\Psi_9(1)=b_1^{\mu(9)}b_3^{\mu(3)}b_9^{\mu(1)}=(-1)^0(-2)^{-1}2^1=-1.\]
Thus
\[\chi(s_1)\chi(s_2)\chi(s_3)=\chi(-1)=\pm 1\quad\hbox{if}\quad p\equiv\pm 1\,{\rm mod}\,(4),\]
so (as in the case $n=7$) the number $l_p$ of outer regular maps is even or odd for $p\equiv\pm 1$ mod~$(4)$.

 If $p=2$ then $d=3$, so $G={\rm PSL}_2(8)$, and there is a unique map R15.1 of extended type $\{3,9\}_{14}$. It is inner regular, since every element of ${\mathbb F}_8$ is a square, with non-orientable quotient N16.1 of type $\{3,9\}_7$.
If $2<p\equiv\pm 2$ or $\pm 4$ mod~$(9)$ the unique map is inner or outer regular as $p\equiv\pm 1$ mod~$(4)$. If $p\equiv\pm 1$ mod~$(9)$ then as in the case $n=7$, if $p\equiv 1$ mod~$(4)$ then $l_p=0$ or $2$, so the number $k_p=3-l_p$ of inner regular maps is $3$ or $1$, whereas if $p\equiv -1$ mod~$(4)$ then $l_p=1$ or $3$, so that $k_p=2$ or $0$. 

For each $p\equiv \pm 1$ mod~$(9)$ there are three maps, with $d=1$.  As will be shown in Section~\ref{sec:Galois} the Galois group of $f_2(s):=f_1(s^2)$ is ${\rm C}_2\wr{\rm C}_3$, and it follows that the relative densities of the sets $\Sigma_k$ are given by the binomial distribition $1/8$, $3/8$, $3/8$ and $1/8$ for $k=0, 1, 2, 3$.

If $p=19$, for instance, then $g=96$, and these maps $\M_j$, with $G={\rm PSL}_2(19)$, are R96.1--3 in the catalogues, with Petrie lengths 20, 18 and 10. Taking $\zeta=4$ as a primitive $9$th root of $1$ we have traces $t_1=\pm 9$ giving $s_1=-2=6^2$, $t_2=\pm 3$ giving $s_2=13$ (a non-square), and $t_4=\pm 7$ giving $s_4=11=7^2$, so two maps are inner regular while one is outer regular. In the catalogues there are two possible non-orientable quotients N97.1 and N97.2 of the required genus $97$ and type $\{3,7\}$, confirming that $k_{19}=2$; they have Petrie lengths 9 and 10, so these are quotients of R96.2 and R96.3 respectively, while R96.1 is outer regular. The fact that $k_{19}=2$ and $l_{19}=1$ is consistent with the above parity result for primes $p\equiv-1$ mod~$(4)$.

If $p=17$ we have the maps R69.1--3, with Petrie lengths 18, 16 and 34. The elements of order $9$ in $G={\rm PSL}_2(17)$ have traces $t=\pm 3$ giving $s=11$ (a non-square), $t=\pm 4$ giving $s=4=2^2$ and $t=\pm 7$ giving $s=5$ (a non-square), so $k_{17}=1$ and $l_{17}=2$. This is consistent with the above parity result for primes $p\equiv 1$ mod~$(4)$. There is only one possible non-orientable quotient of genus $70$, namely N70.1 of Petrie length 17, so this is a quotient of the inner regular map R69.3.

If $p=37$ there are three maps of type $\{3,9\}$  and genus $704$, namely R704.1--3, all with Petrie length $38$. Taking $-3$ as a primitive $9$th root of $1$ we find that the elements of order $9$ in ${\rm PSL}_2(37)$ have traces $t_1=9$ giving $s_1=-4$ (a square), $t_2=5$ giving $s_2=15$ (a non-square) and $t_4=-14$ giving $s_4=-8$ (a non-square). There is one non-orientable map N705.1 of genus $705$ and type $\{3,9\}$, confirming that $k_{37}=1$ and $l_{37}=2$. It is unusual that in this case all three maps R704.1--3 have the same Petrie length, so that it is not clear from~\cite{Pot} which one of them is inner regular; however, a MAGMA computation by Conder~\cite{Con25} of orientable double covers of non-orientable regular maps up to genus 1502 shows that it is R704.3.
 
\medskip

\noindent{\bf Example 8} Let $n=11$, so that $\varphi(n)/2=5$. If $p\equiv\pm 1$ mod~$(11)$ there are five non-isomorphic maps, with $d=1$, whereas if $p\equiv\pm 2, \pm 3, \pm 4$ or $\pm 5$ mod~$(11)$ there are five mutually isomorphic maps, with $d=5$.
 We have
\[\Psi_1(2x)\Psi_{11}(2x)=2(T_6(x)-T_5(x))\]
\[=2((32x^6-48x^4+18x^2-1)-(16x^5-20x^3+5x))\]
\[=(2x-2)(32x^5+16x^4-32x^3-12x^2+6x+1),\]
so that 
\[\Psi_{11}(x)=x^5+x^4-4x^3-3x^2+3x+1,\]
and hence $\Psi_{11}(1)=-1$. Alternatively,
\[\Psi_{11}(1)=b_1^{\mu(11)}b_{11}^{\mu(1)}=(-1)^{-1}1^1=-1.\]
Thus
\[\prod_j\chi(s_j)=\chi(-1)=\pm 1\quad\hbox{as}\quad p\equiv \pm 1\,{\rm mod}\,(4),\]
as in the cases $n=7$ and $9$, so the conclusions about the parities of $k_p$ and $l_p=5-k_p$ are the same as in those cases. The Galois group is now ${\rm C}_2\wr{\rm C}_5$, rather than ${\rm C}_2\wr{\rm C}_3$, but by a similar argument the relative densities of the sets $\Sigma_k$ are
given by the binomial distribution $1/32$, $5/32$, $10/32$, $10/32$, $5/32$ and $1/32$ for $k=0,\ldots, 5$.

If $p=2$ then $d=5$ and there is one orientably regular map; this is R1241.1 in~\cite{Pot}, of extended type $\{3,11\}_{62}$. It is inner regular, with non-orientable quotient N1242.1 of extended type $\{3,13\}_{31}$, so $k_2=1$ and $l_2=0$. Unfortunately, the next case, with $p=23$, $d=1$ and $g=2761$, is too large for the catalogues.

\medskip

\noindent{\bf Example 9} If $n=13$ then $\varphi(n)/2=6$ and
\[\Psi_{13}(x)=x^6+x^5-5x^4-4x^3+6x^2+3x-1,\]
so that $\Psi_{13}(1)=1$. Thus in this case
\[\prod_j\chi(s_j)=\chi(1)=1\]
for all $p$, so that $l_p$ is even.

The primes $p=2$, $3$, $5$ and $53$ have orders $d=6$, $3$, $2$ and $1$ in ${\mathbb Z}_{13}^*/\{\pm 1\}$, so we obtain maps of type $\{3,13\}$  with groups ${\rm PSL}_2(q)$ for $q=2^6$, $3^3$, $5^2$ and $53$. The first and last maps have genera too large to be in any current  catalogues, but when $q=5^2$ we have three orientably regular maps R351.1--3 of type $\{3,13\}$ in~\cite{Pot}, with Petrie lengths $8$, $24$ and $26$; they have a single non-orientable quotient N352.1 of type $\{3,13\}_{13}$,  in~\cite{Pot} so $k_5=1$, with just R351.3 inner regular, and $l_5=2$. When $q=3^3$ we have $s_j=3-t_j^2=-t_j^2$, which is a non-square since $3^3\equiv -1$ mod~$(4)$, so both maps are outer regular; in fact, there are two orientably regular maps R442.2 and R442.3 of types $\{3,13\}_{26}$ and $\{3,13\}_{28}$ in~\cite{Pot}, but they have no non-orientable quotients of genus $443$ listed there; thus $k_3=0$ and $l_3=2$. 

\medskip

\noindent{\bf Example 10} If $n=15$ then $\varphi(n)/2=4$ and
\[\Psi_{15}(x)=x^4-x^3-4x^2+4x+1,\]
so that $\Psi_{15}(1)=1$. Thus in this case we also have
\[\prod_j\chi(s_j)=\chi(1)=1\]
for all $p$, so that $l_p$ is even.

If $p=2$ then $d=4$ and there is a single orientably regular map R205.4 in \cite{Con09, Pot} of genus $205$ and extended type $\{3,13\}_{34}$; since every element of ${\mathbb F}_{2^4}$ is a square, it is inner regular, so $k_2=1$ and $l_2=0$. Its non-orientable quotient is N206.1 of extended type $\{3,9\}_{17}$.

If $p=29$ then $d=1$ and there are four maps of genus $610$ and type $\{3,15\}$ in~\cite{Pot}, namely R610.1--4 with 
Petrie lengths $10, 14, 30, 30$; they have two non-orientable quotients, N611.1 and N611.2, with Petrie lengths 7 and 15, so $k_{29}=l_{29}=2$. Clearly R610.2 is inner regular, with quotient N611.1, while Conder's computation~\cite{Con25} shows that R610.3 (rather than R610.4) is also inner regular, with quotient N611.2.

If $p=31$ then again $d=1$ and there are four maps of genus $745$
 in~\cite{Pot}, namely R745.1--4 with Petrie lengths  8, 16, 30, 32. Taking $\zeta=10$ as a primitive $15$-th root of $1$ we have traces $t_j=7$, $16$, $6$ and $3$ for $j=1$, $2$, $4$ and $8$, giving $s_j=16=4^2$, $-5$, $-2$ (both non-squares) and $-6=5^2$, so $k_{31}=l_{31}=2$. There are two non-orientable quotients N746.1 and N746.2, of Petrie lengths 16 and 8, so R745.2 and R745.1 are their inner regular double covers. 

\medskip

\noindent{\bf Example 11} If $n=17$ then $\varphi(n)/2=8$ and
\[\Psi_{17}(x)=x^8+x^7-7x^6-6x^5+15x^4+10x^3-10x^2-4x+1,\]
so that $\Psi_{17}(1)=1$. Thus in this case we have
\[\prod_j\chi(s_j)=\chi(1)=1\]
for all $p$, so that $l_p$ is again even.

If $p=2$ then $d=4$ and there are two orientably regular maps R221.1 and R221.2 of genus $221$ and types $\{3,17\}_{30}$ and $\{3,13\}_{34}$, with $G={\rm PSL}_2(2^4)$. They are both inner regular since every element of ${\mathbb F}_{2^4}$ is a square; they have non-orientable quotients N222.1 and N222.2, of types $\{3,17\}_{15}$ and $\{3, 17\}_{17}$. Thus $k_2=2$ and $l_2=0$.  
The next example with $n=17$, namely $p=67$, has genus far too large for the catalogues.

\medskip

\noindent{\bf Example 12} If $n=19$ then $\varphi(n)/2=9$ and
\[\Psi_{19}(x)=x^9+x^8-8x^7-7x^6+21x^5+15x^4-20x^3-10x^2+5x+1,\]
so that $\Psi_{19}(1)=-1$. Thus in this case we have
\[\prod_j\chi(s_j)=\chi(-1)=\pm 1\]
as $p\equiv\pm 1$ mod~$(4)$, and $l_p$ is even or odd respectively.

The catalogue~\cite{Pot} covers orientably regular maps of genus $g=2$ to $1501$. One of the largest of our examples included in it is the following. If $p=37$ then $d=1$ and there are nine orientably regular maps R1444.1--9 of type $\{3,19\}$ and genus 1444, with Petrie lengths 18, 18, 36, 38 (five times) and 74. Five of them have non-orientable quotients, N1445.1--5, so $k_{37}=5$ and $l_{37}=4$. These five maps have Petrie lengths 18, 9, 19, 19 and 37, so this shows that R1444.9 is inner regular with quotient N1445.5. The other identifications are not clear from~\cite{Pot}, but Conder's computation~\cite{Con25} shows that R1444.1, R1444.2, R1444.7 and R1444.8 are inner regular with quotients N1445.1, N1445.2, N1445.3 and N1445.4 respectively.


\section{Examples with $n$ even}\label{sec:evenexamples}

There are two related technical modifications which must be made when applying the criterion with $n$ even. Firstly, if $n$ is coprime to $p$, then $G={\rm PSL}_2(q)$ has elements of order $n$ if and only if $q\equiv\pm 1$ mod~$(2n)$, not mod~$(n)$ as when $n$ is odd. In order for $X$ and $Y$ to generate $G$ we need $q$ to be the smallest such power of $p$, so we take $q=p^d$ where $d$ is the multiplicative order of $p$ in ${\mathbb Z}_{2n}/\{\pm 1\}$. 

Secondly, when assigning traces $t_j=\zeta^j+\zeta^{-j}$ to $XY$, we must take $\zeta$ to be a primitive $2n$-th root of $1$ in ${\mathbb F}_q$,: if it is a primitive $n$-th root then $(XY)^{n/2}=-I$, which is the identity element of $G$, so $XY$ has order $n/2$ rather than $n$. (In fact, we have implicitly done this already with the generator $X$ of order $2$, using the $4$th roots of unity $\pm i$ as its eigenvalues.) Nevertheless, there are still $\varphi(2n)/4=\varphi(n)/2$ traces to consider, since in this case four primitive $2n$th roots $\pm\zeta^{\pm 1}$ all give the same traces.

These two changes, effective doubling $n$ when it is even, make specific examples rarer and more difficult to deal with by hand, so we will just briefly summarise a few examples.

\medskip

\noindent{\bf Example 13} If $n=8$ then $\varphi(n)/2=2$. Since $\Psi_8(x)=x^2-2$ we have $\Psi(1)=-1$, so $l_p$ is even or odd as $p\equiv \pm 1$ mod~$(4)$.

If we take $p=17$ then $d=1$ and there are two maps of type $\{3,8\}$ and genus $52$ with $G={\rm PSL}_2(23)$. Taking $\zeta=3$ as a primitive $16$-th root of $1$ we have traces $t_j=9$ and $5$ for $j=1$ and $3$, so that $s_j=3-t_j^2=-7$ and $-5$ (both non-squares), and hence these maps are both outer regular. In the catalogues they are R52.1 and R52.2; there is no non-orientable quotient of genus $53$, confirming that both maps are outer regular. 

If we take $p=31$ then again $d=1$ and there are two maps, now of genus $311$. The elements of order $8$ in $G={\rm PSL}_2(31)$ are non-diagonalisable, having (like all elements of proper $2$-power order)  traces $t$ for which $t^2-4$ is a non-square. Their eigenvalues lie in the quadratic extension field ${\mathbb F}_{31^2}$, but their traces $t$ can be found within ${\mathbb F}_{31}$ by considering the mapping $t\mapsto t^2-2$, induced by squaring group elements; those $t$ requiring two iterations to reach $0$ (the trace of an involution) correspond to elements of order $8$. These $t$, namely $\pm 5$ and $\pm 14$, give $s=3-t^2=9$ (a square) and $-7$ (a non-square), so the corresponding maps are inner and outer regular.
In the catalogues these maps are R311.1 and R311.2, with Petrie lengths $30$ and $32$; there is only one non-orientable quotient, N312.1 with Petrie length $15$, so R311.1 is inner regular with quotient N312.1, while R311.2 is outer regular.

\medskip

\noindent{\bf Example 14} Let $n=10$, so that $\varphi(n)/2=2$. Since $\Psi_{10}(x)=x^2-x-1$ we have $\Psi_{10}(1)=-1$, so $l_p$ is even or odd as $p\equiv \pm 1$ mod~$(4)$. 

 If we take $p=19$ then $d=1$ and there are two maps of type $\{3,10\}$ and genus $115$ with $G={\rm PSL}_2(19)$. The elements of $G$ of order $10$ have traces $t=\pm 6$ and $\pm 8$, with $s=5=9^2$ and $s=-4$, a non-square, so the corresponding maps are inner and outer regular..These maps are R115.1 and R115.2, of Petrie lengths $18$ and $20$. There is only one non-orientable quotient of genus 116, namely N116.1 of Petrie length $9$, so R115.1 is inner regular with quotient N116.1, and R115.2 is outer regular. .

If we take $p=41$ then again $d=1$ and there are two maps, of genus $1149$. Taking $\zeta=2$ as a primitive $20$th root of $1$ we see that elements of order $10$ in $G={\rm PSL}_2(41)$ have traces $t_1=\pm 18$ and $t_3=\pm 3$, with $s_j=-7$ and $-6$, both non-squares. There are no non-orientable quotients of genus $1150$, confirming that $k_{41}=0$ and $l_{41}=2$.

\medskip

\noindent{\bf Example 15} Let $n=12$, so that $\varphi(n)/2=2$. Since $\Psi_{12}(x)=x^2-3$ we have $\Psi_{12}(1)=-2$, so $l_p$ is even for $p\equiv 1$ or $3$ mod~$(8)$, and odd for $p\equiv -1$ or $-3$ mod~$(4)$.

 If we take $p=23$ then $d=1$ and there are two maps of type $\{3,12\}$ and genus 254, with $G={\rm PSL}_2(23)$. They correspond to elements of order 12 with traces $t
=\pm  3$ or $\pm 8$, and hence $s=-6$ (a non-square) or $8\;(=10^2)$, so they are respectively outer and inner regular. These maps are R254.3 and R254.4 in the catalogues, with Petrie lengths 22 and 24; there is just one non-orientable quotient, N255.3, with Petrie length 11, so R254.3 is inner regular with quotient N255.3, while R254.4 is outer regular.

If we take $p=5$ then $d=2$ and there is a single map, R326.2, of type $\{3,12\}_{26}$ with $G={\rm PSL}_2(5^2)$. We can regard ${\mathbb F}_{5^2}$ as ${\mathbb F}_ 5[x]/(x^2-2)$. Since $\zeta:=2x-1$ is a primitive  $24$-th root of $1$ we can take $t=\zeta+\zeta^{-1}=(2x-1)+(x-2)=-2x+2$ as a trace, so that $s=-2x+1$. This is $-\zeta$; since $\zeta$ is not a square, and $5^2\equiv 1$ mod~$(4)$, $s$ is a non-square, so the map is outer regular. To confirm this, there is no non-orientable quotient of genus 327 in the catalogues.

 \medskip

\noindent{\bf Example 16} If $n=14$ there are $\varphi(14)/2=3$ traces to consider. Since $\Psi_{14}(x)=x^3-x^2-2x+1$ we have $\Psi_{14}(1)=-1$, so $l_p$ is even or odd for $p\equiv\pm 1$ mod~$(4)$.

If $p=29$ then $d=1$ and $g=581$; there are three maps $\M_j$, namely  R581.1--3 of type $\{3,14\}$ with Petrie lengths 10, 28 and 58. We can take $\zeta=2$, of order $28$, giving $t_1=-12$ and $s_1=4$, or $t_3=-10$ and $s_3=-10$, or $t_5=13$ and $s_5=8$, with only $s_1$ a square, so only $\M_1$ is inner regular. Since there is only one non-orientable quotient N582.1 of this type, with Petrie length 29, $\M_1$ must be R581.3.
If $p=3$ then $d=3$, $q=27$ and $g=469$, with just one map; this is R469.3 of extended type $\{3,14\}_{28}$. There is no non-orientable quotient of genus $470$, so the map is outer regular.


\section{Galois groups and densities}\label{sec:Galois}

In Section~\ref{sec:density}, in the case $n=7$, we used the Galois group $G_2\cong{\rm C}_2\wr{\rm C}_3$ 
of the polynomial $f_2$ to give the relative densities of the sets $\Sigma_k$ of primes $p$ for which $k_p=k$ for $k=0,\ldots, 3$. For more general $n\ge 7$, just as $\Psi_n(x)$ is the analogue of the polynomial $f_1$, the analogue of $f_2$ here is $\Psi_n(x^2)$.

 The Galois group $G_2$ of $\Psi_n(x^2)$ is a transitive but imprimitive group of degree $\varphi(n)$, permuting the $\varphi(n)/2$ pairs $\pm \alpha_j$ of roots of $\Psi_n(x^2)$ as $G_1={\rm Gal}\,\Psi_n(x)\cong{\mathbb Z}_n^*/\{\pm 1\}$ permutes the roots $\sigma_j=\alpha_j^2$ of $\Psi_n(x)$. It is a subgroup of the largest such imprimitive group, the wreath product $W={\rm C}_2\wr{\mathbb Z}_n^*/\{\pm 1\}$. The kernel of the action of $W$ on pairs is an elementary abelian $2$-group $E$ of rank $\varphi(n)/2$, with quotient group ${\mathbb Z}_n^*/\{\pm 1\}$; the kernel $N$ of the action of $G_2$ on pairs is a normal subgroup of $W$ contained in $E$.
 
 \begin{lemma}\label{lem:G2=W}
 $G_2=W$ for odd $n=7$, $9$ and $11$, and for even $n=8$, $10$, $12$, $14$, $16$ and $18$.
 \end{lemma}

\begin{figure}[h!]
\begin{center}
 \begin{tikzpicture}[scale=0.4, inner sep=0.7mm]

\node (b) at (4,3) [shape=circle, fill=black] {};
\node (c) at (4,-3) [shape=circle, fill=black] {};
\node (d) at (-4.8,1.5) [shape=circle, fill=black] {};
\node (e) at (-4.8,-1.5) [shape=circle, fill=black] {};
\draw [thick] (5,0) arc (0:360:5);
\draw (0,0) to (4.3,2.5) to (4.3,0);
\draw (0,0) to (-4.3,2.5) to (-4.3,0);
\draw (0,0) to (4.3,-2.5) to (4.3,0);
\draw (0,0) to (-4.3,-2.5) to (-4.3,0);
\draw [thick] (-7,0) to (7,0);
\node at (5.5,-0.8) {$1$};
\node at (-6,-0.8) {$-1$};
\node at (0,-0.8) {$0$};
\node at (.5-4,4.5) {$S^1$};
\node at (7.8,0) {$\mathbb R$};
\node at (5,3.2) {$\zeta_n$};
\node at (5,-3.2) {$\zeta_n^{-1}$};-
\node at (-6.5,2) {$\zeta_n^{(n-1)/2}$};
 \node at (-7,-2) {$\zeta_n^{(n+1)/2}$};-
\node at (2,0.5) {$\pi/6$};-
\draw (3,0) arc (0:30:3);
\end{tikzpicture}

\end{center}
\caption{Some primitive $n$-th roots of $1$, for $n=7$, $9$ or $11$.} 
\label{fig:oddn}
\end{figure}
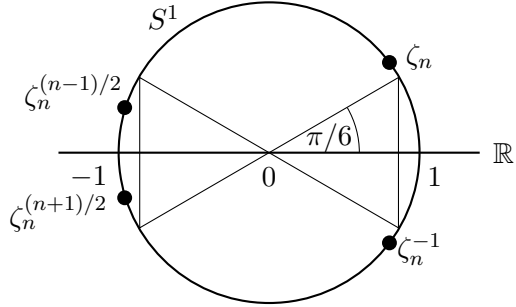

 \noindent{\sl Proof.} As in the case $n=7$, it is sufficient to show that $\Psi_n(x)$ has a single negative root $\sigma_j$, since then $G_2$ contains a permutation $g_j$ transposing the complex conjugate pair of roots $\pm\alpha_j=\sqrt\sigma_j$ of $\Psi_n(x^2)$ and fixing all the other roots, which are real, and the conjugates of $g_j$ generate $E$. Now $\sigma_j=3-\tau_j^2$, so $\sigma_j<0$ if and only if $|\tau_j|>\sqrt 3$. If $n$ is odd this is equivalent to
the real part $\cos(2j\pi /n)$ of the primitive $n$th roots of unity $\zeta^{\pm j}_n$ satisfying $|\cos(2j\pi/n)|>\sqrt 3/2=\cos(\pi/6)$. It is easily seen that for odd $n<13$ only $j=(n-1)/2$ satisfies this, whereas $j=1$ also does for $n\ge 13$ (see Figure~\ref{fig:oddn}).

\medskip

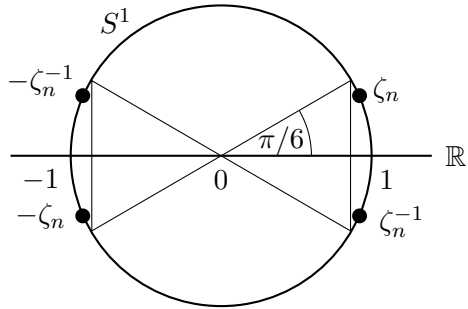
\begin{figure}[h!]
\begin{center}
 \begin{tikzpicture}[scale=0.4, inner sep=0.7mm]

\node (b) at (4.6,2) [shape=circle, fill=black] {};
\node (c) at (4.6,-2) [shape=circle, fill=black] {};
\node (d) at (-4.6,2) [shape=circle, fill=black] {};
\node (e) at (-4.6,-2) [shape=circle, fill=black] {};
\draw [thick] (5,0) arc (0:360:5);
\draw (0,0) to (4.3,2.5) to (4.3,0);
\draw (0,0) to (-4.3,2.5) to (-4.3,0);
\draw (0,0) to (4.3,-2.5) to (4.3,0);
\draw (0,0) to (-4.3,-2.5) to (-4.3,0);
\draw [thick] (-7,0) to (7,0);
\node at (5.5,-0.8) {$1$};
\node at (-6,-0.8) {$-1$};
\node at (0,-0.8) {$0$};
\node at (.5-4,4.5) {$S^1$};
\node at (7.8,0) {$\mathbb R$};
\node at (5.5,2.2) {$\zeta_n$};
\node at (6,-2.2) {$\zeta_n^{-1}$};-
\node at (-6,2.5) {$-\zeta_n^{-1}$};
 \node at (-6,-2) {$-\zeta_n$};-
\node at (2,0.5) {$\pi/6$};-
\draw (3,0) arc (0:30:3);
\end{tikzpicture}

\end{center}
\caption{Some primitive $2n$-th roots of $1$, for even $n>6$.} 
\label{fig:evenn}
\end{figure}

The argument is similar for even $n$, except that we now work with $2n$th roots of unity, and quadruples $\pm\zeta_n^{\pm j}$ of them, instead of pairs. In this case $|\cos(2j\pi/2n)|$ is maximised, and greater than $\cos(\pi/6)$, when $j=1$ (see Figure~\ref{fig:evenn})), while the second largest value is attained by taking the least $j>1$ coprime to $2n$, and this value is less than $\cos(\pi/6)$ for even $n\le 18$. \hfill$\square$

\medskip

The argument used above for $n$ odd breaks down when $n\geq 13$ since then $\cos(2\pi/n)>\cos(\pi/6)$, so that complex conjugation  transposes $\pm\alpha_j$ for $j=1$ and $(n-1)/2$ (and possibly other $j$) rather than just $(n-1)/2$. Similarly the argument for $n$ even breaks down when $n=20$, since taking $j=3$ gives $\cos(2j\pi/2n)>\cos(\pi/6)$. A computer search by Martin Ma\v caj has shown that for odd $n\ge 13$ it is more  usual for $G_2$ to be a proper subgroup of $W$, the first exception being $n=19$, where complex conjugation induces the odd permutation $g_1g_8g_9$.

When the Galois group $G_2$ is $W$, as in Lemma~\ref{lem:G2=W}, one can apply the density theorem as in the case $n=7$ (see Section~\ref{sec:density}) to show that the relative densities of the sets $\Sigma_k$ ($k=0,\ldots, r=\varphi(n)/2$) are  ${r\choose k}/2^r$, forming a binomial distribution. Among such $n$ we have $r=2$ when $n=8$, $10$ or $12$; $r=3$ when $n=7$, $9$, $14$ or $18$; $r=4$ when $n=16$, $r=5$ when $n=11$, and $r=9$ when $n=19$.

Even when $G_2$ is a proper subgroup of $W$, if the subgroup $N<E$ is known one can obtain the relative densities of the sets $\Sigma_k$, since all the elements of $G_2$ with fixed points on the roots lie in $N$. For example, since $N$ is normal in $G_2$ there is only one maximal subgroup $N$ of $E$ which can arise; this is the even subgroup $E^+$ of $E$, consisting of the elements of $E$ inducing even permutations of the roots, so that the relative densities are ${r\choose k}/2^{r-1}$ for even $k$. This situation arises when $n=13$ or $15$, for instance, with $r=6$ and $4$ respectively.


\section{Generalisation to other types}\label{sec:m>3}

One could ask similar questions about the maps or hypermaps of other types besides those considered here. In this case the approach used in Section~\ref{sec:alternative} meets the difficulty that, by a theorem of Takeuchi~\cite{Tak}, triangle groups of only finitely many types are arithmetic groups. However, the approach used in Sections~\ref{sec:f1} and \ref{sec:density} can be transferred from maps of type $\{3,n\}$ to those of type $\{m,n\}$ for any fixed $m>3$. The general principles are the same as for $m=3$, so here we will just outline a short illustrative example.

\medskip

Hall's Criterion has recently been generalised by Conder, Poto\v cnik and \v Sir\'a\v n in~\cite{CPS} to apply to orientably regular hypermaps $\mathcal H$ of arbitrary hyperbolic type $(k,l,m)$ with orientation-preserving automorphism group ${\rm PSL}_2(q)$: the necessary and sufficient condition for the inner regularity (stated equivalently for the existence of a non-orientable regular hypermap) is that
\[-D:=4+\omega_{\kappa}\omega_{\lambda}\omega_{\mu}
-\omega_{\kappa}^2-\omega_{\lambda}^2-\omega_{\mu}^2\]
should be a square in ${\mathbb F}_q$, where $\omega_{\kappa}$, $\omega_{\lambda}$ and $\omega_{\mu}$ are the traces of the standard generators, rotations of $\mathcal H$ through $2\pi/k$, $2\pi/l$ and $2\pi/m$ around a hypervertex, hyperedge and hyperface. (This is the factor on the right-hand side of (\ref{eq:det}), in different notation,) The Macbeath maps of type $\{3,n\}$ can be regarded as hypermaps of type $(n,2,3)$, so that $\omega_{\lambda}=0$, $\omega_{\mu}=\pm 1$ and $-D=3-\omega_{\kappa}^2=3-t^2$ as in the extension of Hall's Criterion.

As an example of this generalisation, let $m=4$, so that we consider maps $\M$ of type $\{4,n\}$. The group $G={\rm PSL}_2(q)$ has elements of this order if and only if $q\equiv\pm 1$ mod~$(8)$, as we will assume throughout this section. Elements of order $4$ have trace $\pm\sqrt 2$, so that $\M$, of type $\{m,n\}$ and with ${\A^+\M}=G$, is inner regular if and only if $s:=2-t^2$ is a square in ${\mathbb F}_q$, where $t$ is the trace of the canonical generator $y$ of order $n$. To have hyperbolic types we require $n\ge 5$, so for simplicity let us take $n=5$, so that $\M$ has genus $1+|G|/40$. We obtain maps $\M$ of type $\{4,5\}$ with ${\rm Aut}^+\M\cong G$ if and only if either $q=p\equiv \pm 1$ or $\pm 9$ mod~$(40)$, with two maps, or $q=p^2$ where $p\equiv \pm 7$ or $\pm 17$ mod~$(40)$, with one map (more precisely, two isomorphic maps). The traces of elements of order $5$ are the roots $t_j=(-1\pm\sqrt 5)/2$ of $f_0(t)=t^2+t-1$, and the corresponding elements $s_j$ are the roots $2-t_j^2=t_j+1=(1\pm\sqrt 5)/2$ of $f_1(s)=s^2-s-1$. Since $s_1s_2=-1$, if $q\equiv 1$ mod~$(4)$ then the two maps $\M_j$ are either both inner or both outer regular, whereas if $q\equiv -1$ mod~$(4)$ then one is inner regular and the other is outer regular.

Since $f_1$ has one positive and one negative root in $\mathbb R$, the Galois group of $f_2(s)=f_1(s^2)$ is ${\rm S}_2\wr{\rm S}_2\cong{\rm D}_4$, with respectively one and two elements fixing either four or two roots, so by imitating the argument in Section~\ref{sec:density} we see that the relative densities of the sets of primes giving $k=0, 1$ or $2$ inner regular maps follow the binomial distribution $1/4, 1/2$ and $1/4$.

\medskip

\noindent{\bf Example 17} First we consider some primes $p\equiv\pm 1$ mod~$(5)$, and hence $p\equiv\pm1$ or $\pm 9$ mod~$(40)$, each yielding two regular maps with $q=p$.

If $p=31$ then $\sqrt 5=\pm 6$, so $t=12$ or $-13$, and hence $s=13$, a non-square mod~$(31)$, or $-12$ ($=9^2$); these give outer and inner regular maps R373.1 and R373.2 respectively, with Petrie lengths $16$ and $32$; there is only one possible non-orientable regular quotient, N374.1 of Petrie length $16$, so this is the quotient of the inner regular map R373.1.

If $p=41$ then $\sqrt 5=\pm 13$, so $t=6$ or $-7$ and hence $s=7$ or $-6$, which are both non-squares. The outer regularity of the corresponding maps, R862.1 and R862.2, is confirmed by the absence from the catalogues of any non-orientable regular maps of type $\{4,5\}$ and genus $863$.

If $p=79$ then $\sqrt 5=\pm 20$, so $s=-29=34^2$ or $30$ (a non-square), giving inner and outer regular maps of genus $6163$, too large for current catalogues.

If $p=89$ then $\sqrt 5=\pm 19$, with $s=10=30^2$ or $-9=13^2$, giving two inner regular maps of genus $8812$.

Next we consider some primes $p\equiv\pm 2$ mod~$(5)$, and hence $p\equiv \pm 7$ or $\pm 17$ mod~$(40)$, each yielding a single regular map with $q=p^2$.

If $p=7$ then $q=49$ and one can check that the elements $s=(1\pm\sqrt 5)/2$ are not squares in ${\mathbb F}_q={\mathbb F}_7(\sqrt 5)$, so the single corresponding map, namely R1471.2, is outer regular. Indeed, the catalogues confirm this by showing that there is no non-orientable regular map of type $\{4,5\}$ and genus $1472$.

However, if $p=17$ then the elements $(1\pm \sqrt 5)/2$ have square roots $7\mp 3\sqrt 5$ in ${\mathbb F}_q$, so the corresponding map of genus $301\,717$ is inner regular.


\section{Other generalisations}\label{sec:general}

In addition to considering the problem raised at the end of Sections~\ref{sec:density} and \ref{sec:alternative}, of characterising the primes in each set $\Sigma_k$, either for $n=7$ or for arbitrary $n$, one can generalise the problem addressed in this paper in other ways.

Building on joint work on maps with Ma\v caj and \v Sir\' a\v n in~\cite{JMS}, a similar approach has recently been used with the same authors in~\cite{JMS25} to show that for each hyperbolic triple there exist infinitely many non-orientable regular hypermaps of that type. (This is significantly harder to prove than the corresponding result for orientable hypermaps, which is an easy consequence of the residual finiteness of triangle groups.) The hypermaps constructed have automorphism group ${\rm PSL}_2(p)$ for sets of primes $p$ of positive density, those where the  inner regularity of an orientably regular hypermap yields a non-orientable hypermap of the same type and with the same group. This approach is both more general than that in the present paper, in that it applies to hypermaps of every hyperbolic type, but also less general in that it applies only in the case $q=p$ and that it does not yield exact or relative densities. Thus the two papers are similar but complementary.

Alternatively, one could consider other families of Hurwitz groups. For example, Conder~\cite{Con80} has proved that all but finitely many alternating groups ${\rm A}_n$ are Hurwitz groups. Many (perhaps most) of the associated Hurwitz maps are chiral, rather than fully regular, but it still makes sense to call the latter inner or outer regular as the full automorphism group is ${\rm A}_n\times{\rm C}_2$ or the symmetric group ${\rm S}_n$ (the only possibilities), and to consider the Klein surfaces appearing in the inner regular cases. For instance, Etayo and Mart\'\i nez~\cite{EM} have shown that, of the five smallest alternating groups which are Hurwitz groups, namely those with $n=15$, $21$, $22$, $28$ and $29$, only ${\rm A}_{15}$ and ${\rm A}_{28}$ also arise as $H^*$-groups. The rapidly increasing size of the alternating groups, together with their lack of any simple representation over a finite field, makes this a much harder problem than that for the Macbeath-Hurwitz groups considered here.


\section{Appendix: computations by Conder}\label{sec:Conder}

Here are the results of Conder's computations with MAPLE mentioned at the end of Section~\ref{sec:regularity}, involving the first $400$ primes $p\equiv\pm 1$ mod~$(7)$. As explained in that section, the sets $\Sigma_k$ consist of these primes $p$ for which there are $k=0$, $1$, $2$ or $3$ inner regular Hurwitz maps associated with ${\rm PSL}_2(p)$. Here each is partitioned into two subsets $\Sigma_k^{\pm}$ as $p\equiv\pm 1$ mod~$(7)$, in order to reveal any possible bias between these two congruence classes; the data here does not show significant evidence of this in any of the sets $\Sigma_k$.

Recall that by Proposition~\ref{prop:k} the primes in $\Sigma_0\cup\Sigma_2$ and $\Sigma_1\cup\Sigma_3$ are congruent to $\pm 1$ mod~$(4)$ respectively.

\bigskip

\noindent $\Sigma_0^+$: 239, 379, 491, 547, 1051, 1583, 2143, 3319, 3823, 3907, 4159, 4271, 4523, 4663, 5503, 5867, 5923, 6427, 6959, 7043, 8443, 9227 (22 primes).

\medskip

\noindent $\Sigma_0^-$: 167, 251, 1399, 1511, 1931, 1987, 2351, 3023, 3331, 3359, 4003, 4283, 4339, 4759, 4871, 5179, 5683, 6803, 7307, 8623, 8707, 9043, 9127, 9239, 9491, 9631 (26 primes).

\bigskip

\noindent $\Sigma_1^+$: 29, 113, 197, 281, 337,  421,  449,  617, 673, 757, 953, 1009,  1429, 1597, 1709, 1877, 1933, 2017, 2129, 2213, 2269, 2297, 2381, 2437,  2633, 2801, 2857, 2969, 3109, 3137, 3221, 3361, 3389, 3529, 3613, 3697, 4201, 4229, 4621, 4649, 4733, 4817, 4957, 5153, 5209, 5237, 5573, 5657, 5741, 5881, 6133, 6217, 6301, 6329, 6469, 6553, 6581, 6637, 7001, 7057, 7309, 7393, 7561, 7589, 8093, 8233, 8317, 8429, 8597, 8681, 8693, 8821, 9157, 9241, 9437, 9521, 9661, 9689, 9829 (79 primes).

\medskip

\noindent $\Sigma_1^-$: 13, 41, 97, 349,  433,  461, 601, 769,  797, 853,  1021, 1049, 1217, 1301, 1609, 1637, 1693, 1721, 1777, 1861, 1889, 1973, 2029, 2113, 2141, 2309, 2477, 2617, 2729, 2953, 2897, 3121, 3541, 3821, 3989, 4073, 4129, 4157, 4409, 4493, 4549, 4801, 4969, 5081, 5333, 5417, 5557, 5641, 5669, 6089, 6173, 6229,  6397, 6733, 6761, 7013, 7069, 7321, 7349,  7433, 7489, 7517, 7573, 7853, 7993, 8161, 8273, 8329, 8861, 9001, 9029, 9281, 9337, 9421, 9533 (75 primes).
 
 \bigskip

\noindent $\Sigma_2^+$: 43, 71, 127, 211, 463, 631, 659, 743, 827, 883, 911, 967, 1163, 1303, 1471, 1499, 1667, 1723, 2003, 2087, 2311, 2339, 2423, 2591, 2647, 2731, 2843, 2927, 3011, 3067, 3347, 3571, 3739, 3767, 3851, 4019, 4243, 4327, 4691, 4831, 4943, 4999, 5167, 5279, 5419, 5531, 5783, 5839, 6007, 6091, 6203, 6287, 6343, 6679, 6763, 6791, 7127, 7211, 7351, 7547, 7603, 7687, 7883, 8191, 8219, 8387, 8527, 8779, 8807, 8863, 9059, 9199, 9283, 9311, 9479, 9619, 9787, 9871 (78 primes).

\medskip

\noindent $\Sigma_2^-$: 83, 138, 223, 307, 419, 503, 587, 643, 727, 811, 839, 1063, 1091, 1231, 1259, 1427, 1483, 1567, 1847, 2099, 2239, 2267, 2659, 2687, 2939, 3079, 3163, 3191, 3499, 3527, 3583, 3779, 3863, 3919, 3947, 4423, 4451, 4507, 4591, 4703, 4787, 5011, 5039, 5347, 5431, 5711, 5851, 5879, 6047, 6131, 6271, 6299, 6551, 6607, 6691, 6719, 6971, 7027, 7559, 7643, 7699, 7727, 7867, 7951, 8147, 8231, 8287, 8539, 8819, 9323, 9463, 9547, 9743 (73 primes). 
 
 \bigskip
 
\noindent $\Sigma_3^+$: 701, 1093, 1289, 1373, 2521, 2549, 2689, 3557, 4397, 4481, 4789, 6833, 6917, 7253, 7477, 7673, 7757, 7481, 8009, 8513, 8737, 8849, 8933, 9857 (24 primes).

\medskip

\noindent $\Sigma_3^-$: 181, 293, 881, 937, 1553, 2281, 2393, 3037, 3373, 3457, 2709, 3793, 3877, 4241, 4297, 5501, 6257, 6481, 7237, 7741, 7937, 8581, 8609 (23 primes).
 
 \bigskip
 

\bigskip

\noindent School of Mathematical Studies\\
University of Southampton\\
Southampton SO17  1BJ\\
UK\\
{\tt G.A.Jones@maths.soton.ac.uk}\\

\end{document}